\documentclass[11pt, oneside]{article}   	
\usepackage{geometry}                		
\usepackage[utf8]{inputenc}

\newcommand\blfootnote[1]{%
  \begingroup
  \renewcommand\thefootnote{}%
  \footnote{#1}%
  \endgroup
}

\usepackage{graphicx}				
\usepackage{amssymb}
\usepackage{amsfonts}
\usepackage{mathrsfs}  
\usepackage{amsbsy}
\usepackage{amsmath}
\usepackage{color}
\usepackage{subfig}

\def\c{{v_\text{net}}}

\title{Deconvolution of inclined channel elutriation data to infer particle size distribution}
\author{Jeffrey A. Hogan${}^{1}$\thanks{corresponding author}, Simon Iveson${}^2$, Jason Mackellar${}^3$, Kevin Galvin${}^4$}
\date{%
     {\tt jeff.hogan@newcastle.edu.au}\\
     {\tt jason.mackellar@newcastle.edu.au}\\
    ${}^{1,3}$ School of Information and Physical Sciences\\
    University of Newcastle\\
    Callaghan, NSW 2308\\
    Australia\\[2ex]
    {\tt  simon.iveson@newcastle.edu.au}\\
    {\tt  kevin.galvin@newcastle.edu.au}\\
    ${}^{2,4}$ ARC Centre of Excellence for Enabling Eco-efficient Beneficiation of Minerals\\
Newcastle Institute for Energy and Resources\\
Department of Chemical Engineering\\
University of Newcastle\\
Callaghan, NSW 2308\\
Australia\\[2ex]%
    \today
}
\begin{document}

\maketitle

\abstract{In this paper we investigate the application of optimisation techniques in the deconvolution of mineral fractionation data obtained from a mathematical model for the operation of a fluidised bed with a set of inclined parallel channels mounted above. The model involved the transport equation with a stochastic source function and a linearly increasing fluidisation rate, with the overflow solids being collected in a finite number of increments (bags). Deconvolution of this data is an ill-posed problem and regularisation is required to provide feasible solutions. Deconvolution with regularisation is applied to a synthetic feed consisting of particles of constant density that vary in size only. It was found that the feed size distribution could be successfully deconvolved from the bag weights, with an accuracy that improved as the rate acceleration of the fluidisation rate was decreased. The deconvolution error only grew linearly with error in the measured bag masses. It was also shown that combining data from two different liquids can improve the accuracy.\blfootnote{This research was supported by the Australian Government through the Australian Research Council's Discovery Projects funding scheme (project DP200102122).}}

\smallskip\noindent
{\bf keywords:\,}{deconvolution, inverse problems, optimisation, regularisation, fractionation}

\smallskip\noindent
{\bf AMS subject classifications:\,}{45Q05, 47A52, 65R30} 

\maketitle

\vskip0.2in
\section{Introduction}
\label{sec: introduction}
The need to transition the global economy to a renewable low-carbon basis is generating a growing demand for critical minerals/metals. The existing high-grade deposits of these minerals are rapidly dwindling, so future demand will need to be met from ores with lower grade and/or a finer mineral grain size. This presents a growing challenge for the minerals processing industry which is driving the search for improved beneficiation methods. This in turn requires effective methods for objectively evaluating how well these new technologies perform.

The concentration of a valuable mineral in a particle often correlates closely with particle density. Hence measuring the density distribution of the feed particles enables one to predict the best possible beneficiation performance, which then provides an objective benchmark against which to assess the actual separation performance of a given technology. The density distribution of particles has traditionally been measured using the sink-float method. This uses multiple baths of different liquids with densities that span the range of interest. However, many of the critical minerals have densities above $5000\,$kg/m${}^3$, whereas the traditional organic heavy liquids used in the sink-float method only cover densities up to around $3300\,$kg/m${}^3$ (diiodo methane) or up to $5000\,$kg/m${}^3$ using Clerici solution (aqueous solution of thallium malonate and formate) \cite{Wills}. Unfortunately many of these liquids also have serious toxicity concerns. Solutions of lithium heteropolytungstates (LST) in water are more benign but, even at elevated temperatures, can only reach densities around $3500\,$kg/m${}^3$ \cite{CCC}.
 
A potential alternative method for measuring particle density distribution is to fluidise the sample beneath a set of closely-spaced parallel inclined channels (Figure \ref{fig: classifier diagram}), and then progressively increase the fluid velocity to elutriate increments of particles of increasing density. The inclined channels create a large effective settling area, the so-called Boycott effect \cite{boycott}. Closely spaced channels promote laminar flow with a parabolic velocity profile which results in particles at the channel wall experiencing a local fluid velocity that is roughly proportional to their size. This scaling reduces the influence of particle size on their behaviour \cite{galvin et al 2009}. Very closely spaced channels also generate high shear rates and the resulting shear-induced lift leads to a very powerful density-based separation \cite{galvin-liu 2011}. The continuous version of this system, known commercially as a REFLUX™ Classifier, is being applied to beneficiate a growing range of minerals such as antimony, chromite, iron ore, mineral sands, tantalum, tin and zircon (\cite{chu et al 2020}, \cite{galvin et al 2020b}, \cite{amariei et al}, \cite{galvin et al 2016}, \cite{iveson et al 2024}, \cite{amoash et al 2025}, \cite{galvin et al 2020a}).

We are interested in using this system in batch mode as a laboratory tool to attempt to measure the density distribution of particles. The system starts with an initial inventory of feed particles, and the fluidisation rate is then increased, progressively elutriating denser and/or larger particles into the overflow, which are collected as a function of time. Early work with elutriating fine $0.038-2.00\,$mm coal with water produced yield-ash curves that closely matched those obtained by the sink-float method \cite{walton et al 2010}. Later work using glycerol solutions extended this range of close agreement up to coal particles as large as $16\,$mm \cite{hunter et al 2014}. Water-based batch fractionation has also been used to accurately measure the density distributions and partition curves for a chromite ore covering a maximum particle density of $\approx 4500\,$kg/m${}^3$ with a nominal top size of $0.300\,$mm \cite{galvin et al 2020c}. More recent batch elutriation work has been performed using LST solution ($\approx 2400\,$kg/m${}^3$) as the fluidising liquid. The greater density of the LST solution increases the buoyancy force on the particles which makes it easier for shear-induced lift to occur. This enhances the sensitivity of the separation to particle density \cite{galvin et al 2020c}. This approach has been applied to a nominal $-2.00\,$mm gold-bearing sulfide ore with maximum particle density $\approx 3800\,$kg/m${}^3$, and showed strong performance \cite{lowes et al 2020}. Elutriation data from a set of feed, product and reject samples have also been used to infer the partition curve (separation performance) of a separator \cite{galvin et al 2018}, \cite{iveson et al 2024}.

The above attempts to infer particle density distributions and partition curves from batch elutriation data did not account for the effect of particle size on elutriation behaviour. It was simply assumed that the material was being elutriated in order of increasing particle density. However, whilst the parabolic velocity profile and shear-induced lift phenomena do dampen the influence of particle size on elutriation behaviour, the effect of particle size is not totally eliminated. Hence this simple approach can only work well for a feed with a narrow particle size range. However, it is hoped that an approach that accounts for the confounding influence of particle size might enable more accurate determination of the particle density distribution, even for widely-sized particle samples. 

The aim of this paper is to consider the simple scenario of how well the size distribution of a mono-density particle sample can be inferred (deconvolved) knowing only the masses of the increments elutriated over time as the fluidisation rate is steadily increased, assuming that one has an accurate model of how particles behave in the system. This scenario not of immediate practical relevance since one could easily directly measure the size distribution of a particle sample using a range of methods, such as sieving or laser diffraction. The long-term plan is to extend this approach to the more useful scenario of whether we can infer the full size by density distribution of a particle sample knowing the mass and size distribution of each elutriation increment. However, the simpler scenario investigated in this paper provides a stepping stone to develop mathematical methods for deconvolution of the data and to investigate how accurately deconvolution works. If this simple mono-density scenario cannot be accurately deconvolved, then there is little hope for deconvolving the more complicated size by density scenario. An additional value in this work is that it also enables us to investigate the potential benefits of using the combined data from two different elutriation experiments e.g., using two liquids with different densities.
 
Section 2 presents the simple “forward” model used for predicting the elutriation mass versus time from a given known feed size distribution for the situation where the fluidising rate is steadily increased. Section 3 presents the mathematical method used to deconvolve such elutriation data to infer the original feed size distribution. Section 4 presents and discusses the results when the model in Section 2 is used to generate an artificial set of elutriation data from an assumed feed size data set, and then the output is deconvolved using the approach in Section 3. No attempt is made to experimentally verify the forward model, and it is acknowledged that there may be significant effects that have not been accounted for. However, hopefully it captures enough of the important hydrodynamics in the system to establish whether this deconvolution approach has any potential, and the results of the analysis may also help to inform about the potential benefits of scenarios such as using two different liquids.

\section{Mathematical modelling of the fractionator}
Figure \ref{fig: classifier diagram} presents a diagram of the batch elutriation system modelled in this paper, which consists of a section of parallel inclined channels mounted on top of a vertical section. In this paper, we model a system where the channels are $\ell  = 1\,$m long inclined at an angle $\theta = 70^\circ$ to the horizontal, with a spacing of $z = 1.8\,$mm. In a batch experiment, a known mass $m$ of feed particles would initially be placed in the vertical section. The fluidisation liquid would then be introduced from beneath and flow upwards through the unit, entraining particles into the overflow where they are then collected using a sieve screen or filter sock. These increments of solids would then be dried and weighed to give a set of increment mass $m_i$ versus time $t$ data. Note that if one were analysing a feed with particles that varied in both size and density, then one might also measure the size distribution (e.g., by sieving or optical light scattering techniques) and/or the average density (e.g., using gas or liquid pycnometry) of the particles in each increment to provide additional information for the deconvolution process. However, this paper is confined to the simple scenario of mono-density feed particles and so only the mass of each increment is measured. If successful, then the approach can be extended in the future to the more generalised case of particles of both varying size and density.

\begin{figure}[h]
\begin{center}
\includegraphics[width=12cm,height=7cm]{classifier-diagram.jpg}
\caption{Schematic of the batch elutriation system showing (a) the introduction of fluidising liquid at the base of the unit and capture of overflow solids using a screen or filter sock and (b) channel close-up showing the parabolic velocity profile with a particle lying on the lower surface of the incline with key variables and dimensions defined.}\label{fig: classifier diagram}
\end{center}
\end{figure}

In the earlier batch fractionation work referred to in the introduction, the fluidisation rate was increased in discrete steps. After each step increase, the flowrate was held constant for a period of time, and the solids elutriated in the overflow during that period were collected as a discrete flow sample. However, these sudden step changes in flowrate were observed to cause an initial surge of particles from the vertical section up into the inclined section, and from there into the overflow. It is speculated that hindered settling effects may be significant during these surges, potentially causing the misplacement of particles into the overflow. Hence, this paper will assume that the flow through the system is increased linearly with time, resulting in a superficial velocity in the channels given by:
\begin{equation}
v(t) = v_0 + \lambda (t – t_0)\label{fluid velocity}
\end{equation}
where $t$ is the time passed since the start of the experiment (T), $v_0$ is the initial fluid velocity (L/T), $t_0$ is the time  at which fluid velocity starts to be increased and $\lambda$ is the velocity scaling factor with dimensions of acceleration (L/T${}^2$) which reflects how quickly the fluidisation rate is being increased.

For flow in a rectangular cross-section duct where the width is much wider than the spacing $z$, the hydraulic diameter (defined as four times the cross-sectional area divided by the wetted perimeter) is equal to $2z$ \cite{fox and mcdonald 1998}. Hence the channel flow Reynold’s number is:
$$\text{\it Re}_{\text c}=\frac{2\rho_{\text F}vz}{\mu}$$
where $\rho_{\text F}$ is the fluid density  and $\mu$ is the fluid dynamic viscosity.
For well-established laminar flow ($\text{\it Re}_{\text c} < 2000$), a parabolic velocity profile is formed in the channel. Following the approach of \cite{galvin et al 2009}, \cite{galvin-liu 2011}, it is assumed that the effective upward fluid velocity $v$ experienced by a particle of radius $s$ resting on the wall is given by the local streamline velocity at a distance $s$ from the wall \cite{bird et al 1976}:
$$v=v(s)=\bar v\frac{6s}{z}\left(1-\frac{s}{z}\right)$$  
where $\bar v$ is the average velocity of the fluid in the channels. 

There are numerous empirical equations in the literature for predicting the terminal free settling velocity $v_t$ of spheres of known size and density in quiescent Newtonian fluids under gravitational acceleration $g$. If a particle is lying on the wall, then it is constrained to move only in the tangential (T) direction. Hence the effective gravitational acceleration acting in that direction is assumed to be $g_{\text T}= g\sin\theta$ \cite{galvin-liu 2011}. So choosing the correlation from \cite{zigrang-sylvester 1981} which is valid for particle terminal free settling Reynolds numbers in the range $\text{\it Re}_{\text{t}} < 10^5$, the terminal settling velocity in the tangential direction $v_{\text{tT}}$, of spheres of radius $s$ and density $\rho_{\text{p}}$ is given by:
\begin{equation}
\text{\it Re}_{\text tT}=\frac{2\rho_{\text F} v_{\text{tT}} s}{\mu}=\left[\left(14.51+\frac{[(\rho_{\text p} -\rho_{\text F} )g\sin(\theta )\rho_{\text F} ]^{0.5} 1.83(2s )^{1.5}}{\mu}\right)^{0.5}-3.81\right]^2.\label{Reynolds}
\end{equation}
Using this equation in this manner assumes that the wall has no effect on particle settling behaviour other than lowering the effective gravitational acceleration i.e., it assumes that there is no friction between the particle and the wall, and that the effect of the wall on the fluid flow field around the particle is negligible. So the particle is assumed to slide rather than roll down the wall, which is unlikely to be the case for spherical particles. It is also noted that Equation (\ref{Reynolds}) takes no account of particle shape, but there are correlations available in the literature that do so if that is required (e.g., \cite{haidler-levenspiel 1989}).

The net velocity of particle movement up the wall is then assumed to be given by the difference between the local upwards fluid velocity and the particle’s downwards tangential settling velocity: 
\begin{equation}
v_{\text{net}}= v(s) – v_{\text{tT}}.\label{net velocity}
\end{equation}

If the vertical section beneath the inclined channels is approximated as being well-mixed, then it is reasonable to assume that particles are pushed up into the base of the inclined section at a rate that is proportional to their concentration in the vertical section. It is further assumed here that the rate of particle entry into the inclined channels is also proportional to their net velocity up the channels. Hence  particles of size $s$ and density $\rho_{\text p}$ are assumed to not enter the channels until their net velocity is positive $(v_{\text{net}} > 0)$ and consequently the rate of change of mass $M(s,t)$ of  particles of size $s$ remaining in the vertical section as a function of time for $v_{\text{net}}>0$ is assumed to be given by:
\begin{equation}
\frac{\partial M(s,t)}{\partial t}=-k c(s,t) M(s,t)\label{pde}
\end{equation}
where $c(s,t)=\max\{\c (s,t),0\}$ and
$k$ is the elutriation rate constant for particles which has dimensions inverse length (L${}^{-1}$). Note that if $t_*(s)$ is the critical time at which $v_{\text{net}}$ first becomes positive for particles of size $s$, then 
$$c(s,t)=\begin{cases}
\c (s,t)&\text{ if }t>t_*(s)\\ 
0&\text{ otherwise.}
\end{cases}$$ 
This model is very simple compared to the complex behaviour that actually occurs at the interface between the vertical section and inclined channels. Particles may in reality enter at any point across the channel spacing $z$. As these particles settle downwards onto the lower channel surface, they will pass through regions of varying local velocity and will land at some finite distance up the inclined channel.  Hence there would not in reality be a one-to-one mapping between a particle’s elutriation exit time and the time it entered at the channel base. However, it is here assumed that when particles which satisfy the constraint $v_{\text{net}}>0$ enter the inclined section, they immediately find themselves on the lower wall of the channel at position $x = 0$ and begin to move up the channel at a velocity $v_{\text{net}}>0$ given by (\ref{net velocity}) which increases linearly with time. 

In reality there would be a finite entry length region before the velocity flow profile attains the fully-established laminar parabolic flow profile. Many particles may enter the inclined section well before $v_{\text{net}}>0$, but will then settle onto the channel’s lower surface and slide back down into the vertical section. The inventory of these recirculating particles in the channels is here assumed to be negligible, but in reality particles with $v_{\text{net}}$ below but close to $0$ could build up in significant numbers and create hindered settling effects near the base of the channels. Particles sliding back down could also create an effective negative slip velocity at the channel wall close to the base of the inclines. Further up the inclines, most particles present are moving upwards, which might create an effective upwards slip velocity at the wall. However, in the simple model used here, hindered settling effects have been neglected and particle behaviour is assumed to be independent of the presence of other particles. Hence the elutriation outputs of the different species are summative.
 
Those particles which have just slid down a channel back into the vertical section may also be more likely to then be re-entrained again into the upwards flow back up into the inclined section. Hence the elutriation rate constant $k$ may in reality have some dependence on particle size and density. However, for the purposes of this paper it is taken as being constant for all particles. Shear-induced inertial lift effects have also been assumed to be negligible, so the full model of \cite{galvin-liu 2011} has not been applied.

The validity of many of the assumptions implicit in this model are unknown. However, the aim here is to investigate whether, given a model that describes the system’s elutriation behaviour, it is possible to accurately infer (deconvolve)  the size distribution of the original feed from the output elutriation data.

\subsection{Transport equation}
Let $u=u(x,s,t)$ be the mass distribution of particles of size $s$ at time $t$  along the inclined channels. Here $x$ is the distance along the inclined channels (with $x=0$ the interface between the channels and the fluidised bed), $s$ is the radius of a particle (assumed to be spherical) and $t\geq 0$ is the time elapsed since the system was started. We model the motion of particles along the inclined channels through the transport problem \cite{strauss}:
\begin{align}
\frac{\partial u}{\partial t}+c(s,t)\frac{\partial u}{\partial x}=0&\qquad (s,x,t>0)\notag\\
u(x,s,0)=g(x,s)&\qquad (x,s>0)\label{transport problem}\\
u(0,s,t)=f(s,t)&\qquad (s,t>0).\notag
\end{align}
Here $c(s,t)$ is the (non-negative) velocity function of a particle of size $s$ in the channels at time $t$. The function $g=g(x,s)$ represents the initial distribution of particles along the inclined channels at time $0$ and $f=f(s,t)$ is the source function which describes the way  particles enter the channels from the fluidised bed. Let
$$d_s(t)=\int_0^tc(s,r)\, dr$$
be the distance travelled by a particle of size $s$ over the time interval $[0,t]$. We assume that $c(s,t)\geq 0$ and that for each size $s$, there is a time $t_*(s)$ for which $c(s,t)>0\iff t>t_*(s)$. Under these assumptions, for each fixed $s$, $d_s(t)$ is an increasing function of $t$ and $d_s:[t_*(s),\infty )\to [0,\infty )$ is invertible with inverse $d_s^{-1}:[0,\infty )\to [t_*(s),\infty )$. The solution of the transport problem above is
\begin{equation}
u(x,s,t)=\begin{cases}g(x-d_s(t),s)&\text{ if }x\geq d_s(t)\\
f(s,d_s^{-1}(d_s(t)-x))& \text{ if }x<d_s(t).
\end{cases}\label{u soln}
\end{equation}
Let $M=M(s,t)$ be the mass distribution by size of particles in the fluidised bed at time $t$. The aim is to determine $M(s,0)$ from a finite collection of discrete flow samples. The model for the transfer or particles from the fluidised bed into the inclined plates is given by (\ref{pde}).
The solution of this first-order linear equation is
$$M(s,t)=M(s,0)e^{-kd_s(t)}.$$
Fix $s>0$ and let $m_r(s,t)$ be the mass of particles of size between $s$ and $s+\Delta s$ remaining in the fluidised bed at time $t$, i.e.,
$$m_r(s,t)=\int_s^{s+\Delta s}M(r,t)\, dr.$$
For the same value of $s$, let $m_c(s,t)$ be the mass of particles of size between $s$ and $\Delta s$  in the inclined channels (assumed for the moment to be of infinite length) at time $t$, i.e., $$m_c(s,t)=\int_s^{s+\Delta s}\int_0^\infty u(x,r,t)\, dx\, dr.$$
Since $m_r(s,t)+m_c(s,t)$ is constant, we have
\begin{align*}
0=\frac{\partial m_r}{\partial t}(s,t)+\frac{\partial m_c}{\partial t}(s,t)&=\int_s^{s+\Delta s}\frac{\partial M}{\partial t}(r,t)\,dr+\int_s^{s+\Delta s}\int_0^\infty\frac{\partial u}{\partial t}(x,r,t)\, dx\, dr\\
&=-k\int_s^{s+\Delta s}c(r,t)M(r,t)\, dr-\int_s^{s+\Delta s}\int_0^\infty c(r,t)\frac{\partial u}{\partial x}(x,r,t)\, dx\, dr\\
&=- k\int_s^{s+\Delta s}c(r,t)M(r,t)\, dr+\int_s^{s+\Delta s}c(r,t)u(0,r,t)\, dr\\
&=-k\int_s^{s+\Delta s}c(r,t)M(r,t)\, dr+\int_s^{s+\Delta s}c(r,t)f(r,t)\, dr .
\end{align*}
Since $\Delta s$ is arbitrary, we conclude that 
$kc(s,t)M(s,t)=c(s,t)f(s,t)$
or equivalently,
\begin{equation}
f(s,t){\mathbf 1}_{c(s,t)>0}=kM(s,t){\mathbf 1}_{c(s,t)>0}=kM(s,0)e^{-kd_s(t)}{\mathbf 1}_{c(s,t)>0}\label{f and M}\end{equation}
where ${\mathbf 1}_{c(s,t)>0}$ is the characteristic function
$${\mathbf 1}_{c(s,t)>0}=
\begin{cases}
1&\text{ if }c(s,t)>0\\
0&\text{ else.}
\end{cases}$$

\subsection{An integral equation model for the Fractionator -- the forward problem}

Assume now that the channels have finite length $\ell$ and that a particle is elutriated from the system if it has travelled a distance $\ell$ metres from the bottom of the channels. Assume also that at time $t=0$ there are no particles in the inclined channels, i.e., in (\ref{transport problem}) we have $g\equiv 0$. Let $m_e(t)$ be the total mass that has elutriated in the time period $[0,t]$. Then $m_e(t)=\int_0^\infty\int_\ell^\infty u(x,s,t)\, dx\, ds$ and the rate at which mass is elutriating at time $t>0$ is 
\begin{align*}
\frac{dm_e}{dt}&=\int_0^\infty\int_\ell^\infty\frac{\partial u}{\partial t}(x,s,t)\, dx\, ds\\
&=-\int_0^\infty\int_\ell^\infty c(s,t)\frac{\partial u}{\partial x}(x,s,t)\, dx\, ds=\int_0^\infty c(s,t)u(\ell ,s,t)\, ds.
\end{align*}
Therefore, the mass that elutriates from the system over the time interval $[t_{i-1},t_{i}]$ $(i\geq 1)$ is
\begin{align}
m_{e,i}=m_e(t_{i})-m_e(t_{i-1})&=\int_{t_{i-1}}^{t_{i}}\frac{dm_e}{dt}(t)\, dt\notag\\
&=\int_{t_{i-1}}^{t_{i}}\int_0^\infty c(s,t)u(\ell ,s,t)\, ds\, dt\notag\\
&=\int_{t_{i-1}}^{t_i}\int_0^\infty c(s,t)f(s,d_s^{-1}(d_s(t)-\ell )){\mathbf 1}_{d_s(t)>\ell}\, ds\, dt\notag\\
&=\int_{t_{i-1}}^{t_{i}}\int_0^\infty c(s,t)kM(s,0)e^{-k(d_s(t)-\ell )}{\mathbf 1}_{d_s(t)>\ell}\, ds\, dt\notag\\
&=k\int_0^\infty M(s,0)\int_{t_{i-1}}^{t_{i}}c(s,t)e^{-k(d_s(t)-\ell )}{\mathbf 1}_{d_s(t)>\ell}\, dt\, ds\notag\\
&=\int_0^\infty M(s,0)L_i(s)\, ds\label{int eqn}
\end{align}
where $L_i(s)=k\int_{t_{i-1}}^{t_{i}}c(s,t)e^{-k(d_s(t)-\ell )}{\mathbf 1}_{d_s(t)>\ell}\, dt$ and we have used (\ref{u soln}) and (\ref{f and M}). The integral equation (\ref{int eqn}) relates the unknown mass distribution across size in the fluidised bed at time $t=0$ with the mass of particles that have elutriated between times $t_{i-1}$ and $t_{i}$. It is clear that if $t_{i}<d_s^{-1}(\ell )$, then $L_i(s)=0$. If $t_{i-1}>d_s^{-1}(\ell )$, then 
$$L_i(s)=k\int_{t_{i-1}}^{t_{i}}c(s,t)e^{-k(d_s(t)-\ell )}\, dt=ke^{kl}\int_{d_s(t_{i-1})}^{d_s(t_{i})}e^{-kr}\, dr=e^{k\ell}[e^{-kd_s(t_{i-1})}-e^{-kd_s(t_{i})}].$$
Finally, if $t_{i-1}<d_s^{-1}(\ell )<t_{i}$ then 
\begin{align*}
L_i(s)&=k\int_{t_{i-1}}^{t_{i}}c(s,t)e^{-k(d_s(t)-\ell )}{\mathbf 1}_{d_s(t)>\ell}\, dt\\
&=ke^{k\ell}\int_{d_s^{-1}(\ell )}^{t_{i}}c(s,t)e^{-kd_s(t)}\, dt=ke^{k\ell}\int_\ell^{d_s(t_{i})}e^{-kr}\, dr=1-e^{k(\ell -d_s(t_{i}))}.
\end{align*}
In summary:
\begin{equation}
L_i(s)=\begin{cases}
0&\text{ if }\ell >d_s(t_{i})\\
1-e^{k(\ell -d_s(t_{i}))}&\text{ if } d_s(t_{i-1})\leq\ell\leq d_s(t_{i})\\
e^{k\ell}[e^{-kd_s(t_{i-1})}-e^{-kd_s(t_{i})}]&\text{ if }\ell<d_s(t_{i-1}).
\end{cases}\label{kernel def}
\end{equation}

\subsection{Deconvolution -- the inverse problem}\label{sec: deconvolution}
\label{sec: deconvolution -- the inverse problem}
We now assume that the unknown feed $M(s,0)$ has  the reconstructed form 
\begin{equation}
\bar M(s,0)=\sum_{j=1}^Na_j\phi_j(s)\label{M exp}
\end{equation}
for known functions $\{\phi_j\}_{j=1}^N$ taking values with units of mass per unit length and unknown dimensionless coefficients $\{a_j\}_{j=1}^N$. Under the assumption  (\ref{M exp}),  we have
\begin{equation}
m_e(t_{i})-m_e(t_{i-1})=\int_0^\infty M(s,0)L_i(s)\, ds=\sum_{j=1}^Na_j\int_0^\infty L_i(s)\phi_j(s)\, ds\quad (1\leq i\leq p).\label{disc eqn}
\end{equation}
Given data $m_{e,i}=m_e(t_{i})-m_e(t_{i-1})$ $(1\leq i\leq p)$, let $B$ be the $p\times N$ {\it cross-grammian} matrix with $(i,j)$-th entry
\begin{equation}
B_{ij}=\int_0^\infty L_i(s)\phi_j(s)\, ds.\label{B entries}
\end{equation}
Then equation (\ref{disc eqn}) may be written as
\begin{equation}
B{\mathbf a}={\mathbf m}_e\label{matrix eqn}
\end{equation}
where ${\mathbf m}_e\in{\mathbb R}^p$ has $i$-th entry $m_{e,i}$ and ${\mathbf a}\in{\mathbb R}^N$ has $j$-th entry $a_j$. Integrating equation (\ref{M exp}) with respect to $s$ over the interval $(0,\infty )$ gives (with $m=\int_0^\infty M(s,0)\,ds$ the total mass of the sample and $c_j=\int_0^\infty \phi_j(s)\, ds$ having dimensions of mass)
\begin{equation}
\sum_{j=1}^Nc_ja_j=  m,\label{mass con}
\end{equation}
 another linear condition on ${\mathbf a}$. The deconvolution of the fractionation data $m_e$ amounts to the solution of  (\ref{matrix eqn}) for the coefficient vector ${\mathbf a}$ subject to the constraint (\ref{mass con}).

Typically, $p<N$ so that $B$ is a ``short-fat'' matrix. If
equation (\ref{matrix eqn}) has a solution then it will have infinitely many. We instead solve the quadratic optimisation problem
$$\text{minimise\ }\|B{\mathbf a}-{\mathbf m}_e\|^2+\alpha\|Q{\mathbf a}\|^2\text{ subject to }{\mathbf c}^{\text T}{\mathbf a}= m$$
where ${\mathbf c}\in{\mathbb R}^N$ has $j$-th entry $c_j$ (as in (\ref{mass con})), ${\mathbf c}^{\text T}$ is the transpose of the row vector ${\mathbf c}$,  $Q$ is a suitable $N\times N$ {\it regularisation} matrix and $\alpha >0$ is a (small) regularisation parameter \cite{mueller-siltanen}. Use of the $N\times N$ identity matrix as a regularisation matrix promotes solutions of small norm. Regularisation matrices of the form 
\begin{equation}
Q_{ij}=\begin{cases}1&\text{ if }i=j\\
-1&\text{ if }i=j-1\\
0&\text{ else}
\end{cases}\label{regulariser}
\end{equation}
promote smoothness of the solution \cite{mueller-siltanen}, and it is this matrix that we will use in our calculations.

Even when the functions $\phi_j$ take positive values only, the solutions ${\mathbf a}^*$ of these problems may have negative entries, in which case the estimated distribution $\bar M(s,0)=\sum_{j=1}^Na_j\phi_j(s)$ may take negative values, unlike the actual initial distribution $M(s,0)$. 
In the case where the functions $\phi_j$ are piecewise constant or piecewise linear splines (see Section \ref{sec: inversion}), the condition $a_j\geq 0$ is equivalent to functions of the form (\ref{M exp}) being non-negative, so we impose the inequality constraints $a_j\geq 0$ $(1\leq j\leq N)$ to avoid  solutions that take negative values. The optimisation problem becomes 
\begin{equation}
\text{minimise\ }\|B{\mathbf a}-{\mathbf m}_e\|^2+\alpha\|Q{\mathbf a}\|^2\text{ subject to }{\mathbf c}^{\text T}{\mathbf a}= m\text{ and }a_j\geq 0\text{ for }1\leq j\leq N\label{constrained opt problem}
\end{equation}
where $\alpha >0$ is a (small) regularisation parameter.

\subsection{Model for fluid velocity}
We assume fluid flow in the channels of the form (\ref{fluid velocity}), 
i.e., average fluid velocity  across each channel is constant at $v_0$ from time $t=0$ until time $t=t_0$, after which it increases linearly at the rate $\lambda$. The net velocity of a particle in the channels is of the form 
\begin{equation}
v_{\text{net}}(s,t)=C_1(s)v(t)-C_2(s)\label{net velocity}
\end{equation}
where 
$$C_1(s)=\frac{6s}{z}\left(1-\frac{s}{z}\right);\quad C_2(s)=\text{\it Re}_{\text{tT}}\frac{\mu}{2\rho_{\text F} s}$$
with $\text{\it Re}_{\text{tT}}$ as in (\ref{Reynolds}).
Recall that for each size $s$ we let $t_*(s)$ be the time at which $c(s,t)$ first becomes positive. Since  $c(s,t)$ is a non-decreasing function of $t$ for each fixed $s$, $t_*(s)$ is well-defined. Note that if $C_1(s)v_0-C_2(s)\geq 0$ then $t_*(s)=0$. Suppose now that $C_1(s)v_0-C_2(s)<0$. Then $t_*(s)\geq t_0$ and in fact we have
$$C_1(s)(v_0+\lambda (t_*(s)-t_0))-C_s(s)=0$$
which has solution $t_*(s)=t_0+\frac{C_2(s)-v_0C_1(s)}{\lambda C_1(s)}$. In summary,
$$t_*(s)=\begin{cases} 0&\text{ if }C_1(s)v_0-C_2(s)>0\\
t_0+\frac{C_2(s)-v_0C_1(s)}{\lambda C_1(s)}&\text{ if }C_1(s)v_0-C_2(s)<0.
\end{cases}$$
With this definition of $t_*(s)$, we have
$$c(s,t)=[C_1(s)v_0-C_2(s)]{\mathbf 1}_{t>t_*(s)}.$$
Then the distance travelled  $d_s(t)$  by a particle of size $s$ over the time interval $[0,t]$ is given by 
$$d_s(t)=\int_0^tc(s,r)\, dr=\begin{cases}0&\text{ if }t<t_*(s)\\
\int_{t_*(s)}^tc(s,r)\, dr&\text{ if }t>t_*(s).
\end{cases}$$
We need more precise information about the function $d_s(t)$. The initial net velocity of a particle of size $s$ in the channels is $D(s)=C_1(s)v_0-C_2(s)$ and suppose that $D(s)>0$ so that $t_*(s)=0$. If $0\leq t\leq t_0$, then 
$d_s(t)=\int_0^tD(s)\, dr=tD(s)$.
On the other hand, if $t>t_0$ then 
$$d_s(t)=\int_0^tc(s,r)\, dr=\int_0^{t_0}D(s)\, dr+\int_{t_0}^t[D(s)+\lambda C_1(s)(t-t_0)]\, dr=D(s)t+\frac{\lambda}{2}C_1(s)(t-t_0)^2.$$
In summary, if $D(s)>0$, then 
\begin{align}
d_s(t)&=\begin{cases}D(s)t&\text{ if } t<t_0\\
D(s)t+\frac{\lambda}{2}C_1(s)(t-t_0)^2&\text{ if } t\geq t_0.
\end{cases}\notag\\
&=D(s)t+\frac{\lambda}{2}C_1(s)(\max (t-t_0,0))^2.\label{first case d_s}
\end{align}
Suppose now that $D(s)<0$. Then $t_*(s)=t_0-\dfrac{D(s)}{\lambda C_1(s)}$ and if $t>t_*(s)$ then 
\begin{align*}
d_s(t)&=\int_{t_*(s)}^t[D(s)+\lambda C_1(s)(r-t_0)]\, dr\\
&=D(s)(t-t_*(s))+\frac{\lambda}{2}C_1(s)[t^2-t_*(s)^2-2t_0(t-t_*(s))]\\
&=(t-t_*(s))[D(s)+\frac{\lambda}{2}C_1(s)(t+t_*(s)-2t_0)]=\frac{\lambda}{2}C_1(s)(t-t_*(s))^2.
\end{align*}
In summary, if $D(s)<0$, then 
\begin{align}
d_s(t)&=\begin{cases} 0&\text{ if }t<t_*(s)=t_0-\frac{D(s)}{\lambda C_1(s)}\\
\frac{\lambda}{2}C_1(s)(t-t_*(s))^2&\text{ if }t\geq t_*(s).
\end{cases}\notag\\
&=\frac{\lambda}{2}C_1(s)(\max\{t-t_*(s),0\})^2\label{second case d_s}
\end{align}
We conclude that for all $t>0$, 
\begin{equation}
d_s(t)=\begin{cases}
D(s)t+\frac{\lambda}{2}C_1(s)(\max\{t-t_0,0\})^2&\text{ if }D(s)>0\\
\frac{\lambda}{2}C_1(s)(\max\{t-t_*(s),0\})^2&\text{ if } D(s)\leq 0.
\end{cases}\label{final dist}
\end{equation}
In both cases $D(s)>0$ and $D(s)\leq 0$, the functions $d_s$ defined by (\ref{first case d_s}) and (\ref{second case d_s}) represent increasing functions on $[t_*(s),\infty )$ and therefore have inverses $d_s^{-1}:[0,\infty )\to [t_*(s),\infty )$.

\subsection{Size distribution data}\label{subsec: size dist data}
Given a real particle feed, size distribution data can be obtained  from physical sieving or laser sizing. For each of these methods there are values $\beta >0$, $\alpha >1$ (different for each method) for which the data $\{\gamma_j\}_{j=1}^N\subset [0,1]$ takes the form
$$\text{proportion of the mass of the sample of sizes $s$ in the range }[\beta\alpha^{j}, \beta\alpha^{j+1}]=\gamma_j,$$
which we write as
$$\int_{\beta\alpha^j}^{\beta\alpha^{j+1}}M(s,0)\, ds=\gamma_jm$$
where $m=\int_0^\infty M(s,0)\, ds$ is the total mass of the sample. Hence, a piecewise constant model for the feed is
$$M(s,0)=\frac{m\gamma_0}{\beta}{\mathbf 1}_{[0,\beta ]}(s)+\frac{\bar m}{\beta (\alpha -1)}\sum_{j=1}^N \gamma_j\alpha^{-j}{\mathbf 1}_{[\beta\alpha^{j-1},\beta\alpha^{j}]}(s)$$
where $\gamma_0=\frac{1}{m}\int_0^\beta M(s,0)\, ds$ and $\beta\alpha^N$ is the upper bound on the size of particles in the sample. For the purposes of the computations we report on in this paper, we use parameters  $\alpha=21.2\times 10^{-6}\,$m for the lower size limit and $\beta =1.1369$ for the scaling of the size increments.

In this paper we consider synthetic feeds only. In particular, given the nature of the grinding processes used on mineral samples to prepare them for analysis, and in an effort to produce a realistic feed, we assume that our feed is a Rosin-Rammler distribution \cite{rosin-rammler}, namely
\begin{equation}
M(s)=A\left(\frac{s}{z}\right)e^{\log(0.2)(2s/(zP))^{2}}\label{RR}
\end{equation}
where $s$ is particle radius, $z$ is the width of the channels, $P=0.125$ and $A$ is chosen so that $\int_0^\infty M(s)\, ds=1$. A plot of this distribution is given in Figure \ref{fig: 5 bag reconstructions runtime =2.67 hours}.

\section{Deconvolution from synthetic data}\label{sec: inversion}

In this section we demonstrate the process of deconvolving (synthetic) bag mass data to produce approximations of the initial size distribution of the sample. 

Published elutriation experiments in inclined channels have used water and aqueous LST solutions \cite{walton et al 2010}, \cite{lowes et al 2020}, so we have chosen typical values for the density and viscosity of these two liquids to use in our simulations, as shown in Table \ref{table: viscosity and density}.

\begin{table}[h]\begin{center}
\begin{tabular}{|c|c|c|}
\hline
&fluid viscosity $\mu$ (Pa\, s)&fluid density $\rho_{\text F}$ (kg/m${}^3$)\\
\hline
water&0.001&1000\\
LST solution&0.0035&2200\\
\hline
\end{tabular}
\caption{Viscosity and density parameters for water and LST solutions used in the numerical experiments.}\label{table: viscosity and density}
\end{center}
\end{table}

\subsection{Numerics of the forward problem}\label{subsec: forward prob}

We assume that the particle density in our experiment is constant, and equal to that of silica, i.e., $\rho_{\text p}=2650\,$kg/m${}^3$ and that the  rate constant is $k=0.0138\,$m${}^{-1}$. We also assume $t_0=0$, i.e., the linear increase in the flow rate starts immediately, that the initial average fluid velocity at time $t=0$ is $v_0=0\,$m/s. Based on typical geometries used in the literature, we assume $z=1.8\times 10^{-3}\,$ m and the length of the inclined channels is $\ell =1\,$m. To generate bag mass data for the forward problem, we assume a feed distribution of the form  (\ref{RR}).

\begin{figure}%
    \centering
    \subfloat{{\includegraphics[width=8cm]{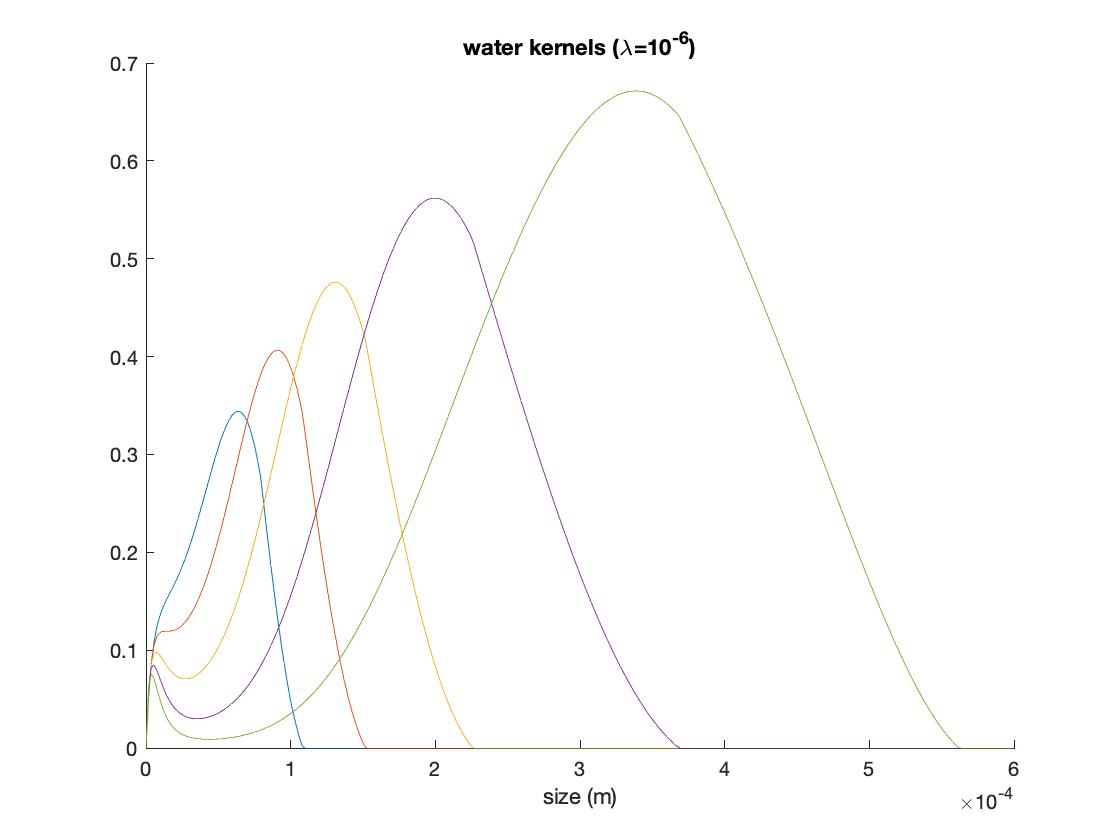} }}%
    \subfloat{{\includegraphics[width=8cm]{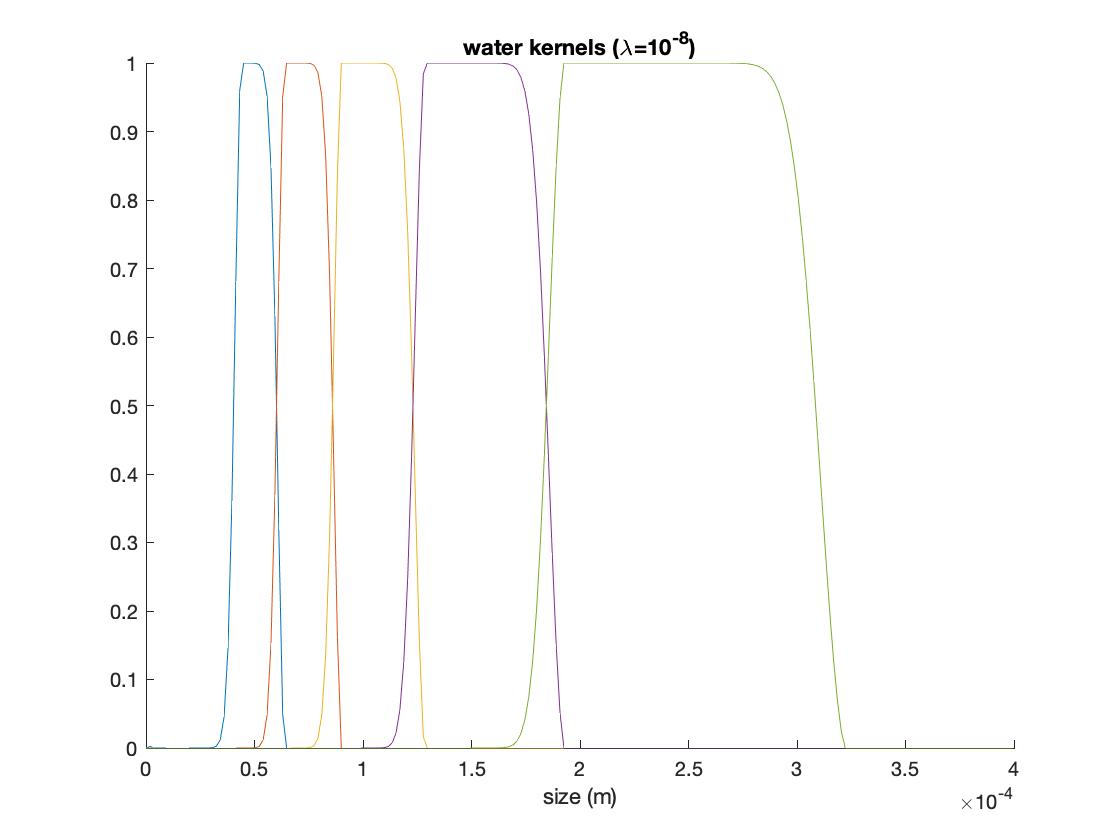} }}%
    \caption{Plots of $5$ water kernel functions $L_i^1$ associated with times evenly spaced between $T_0$ (the time at which $5\%$ of the mass has elutriated) and $T_4$ (the time at which $95\%$ of the mass has elutriated). The plot on the left shows kernels with $\lambda =10^{-6}\,$m/s${}^2$ and the plot on the right shows kernels with $\lambda =10^{-8}\,$m/s${}^2$}%
    \label{fig: kernel plots}%
\end{figure}

With these parameters, we compute the kernels $L_i^1$ associated with water and $L_i^2$ associated with the LST solution as in (\ref{kernel def}). For the purpose of visualising these kernels, we let $T_0$ be the time at which $5\%$ of the feed mass has elutriated and $T_N$ be the time at which $95\%$ of the feed mass has elutriated. We take $p$ bags at uniformly spaced times between $T_0$ and $T_N$.  The mass we expect in the $i$-th bag is given by 
\begin{equation}
m_{e,i}^\varepsilon=\int_0^\infty M(s,0)L_i^\varepsilon(s)\, ds .\label{model masses}
\end{equation}
with $M=M(s,0)$ the size distribution of the feed and $L_i^\varepsilon$ $(i\geq 1,\ \varepsilon\in\{1,2\})$ the kernel function associated with the $i$-th bag and the relevant fluid ($\varepsilon=1$ for water, $\varepsilon =2$ for the LST solution). The first $5$ of these kernels associated with water runs of the experiment are plotted in Figure \ref{fig: kernel plots} (with $\lambda =10^{-6}\,$m/s${}^2$ on the left and $\lambda =10^{-8}\,$m/s${}^2$ on the right). These are indexed by increasing intersections with the $s$ axis.  For both solutions, the small-time kernels $L_i^\varepsilon$ are concentrated on small sizes, meaning that we expect only small particles to be elutriated early in the experiment. As time goes by, the support of the kernels broaden so that larger particles eventually elutriate. The nature of the kernels vary greatly with the value of $\lambda$. For high values of $\lambda$, the supports of the  kernels (see Figure \ref{fig: kernel plots}) have large overlaps so that each bag will be composed of particles with a broad range of sizes. For low values of $\lambda$, the supports of the kernels have much less overlap so that each bag is composed of particles with a narrow range of sizes.

\begin{figure}[h]
\begin{center}
\includegraphics[width=12cm,height=10cm]{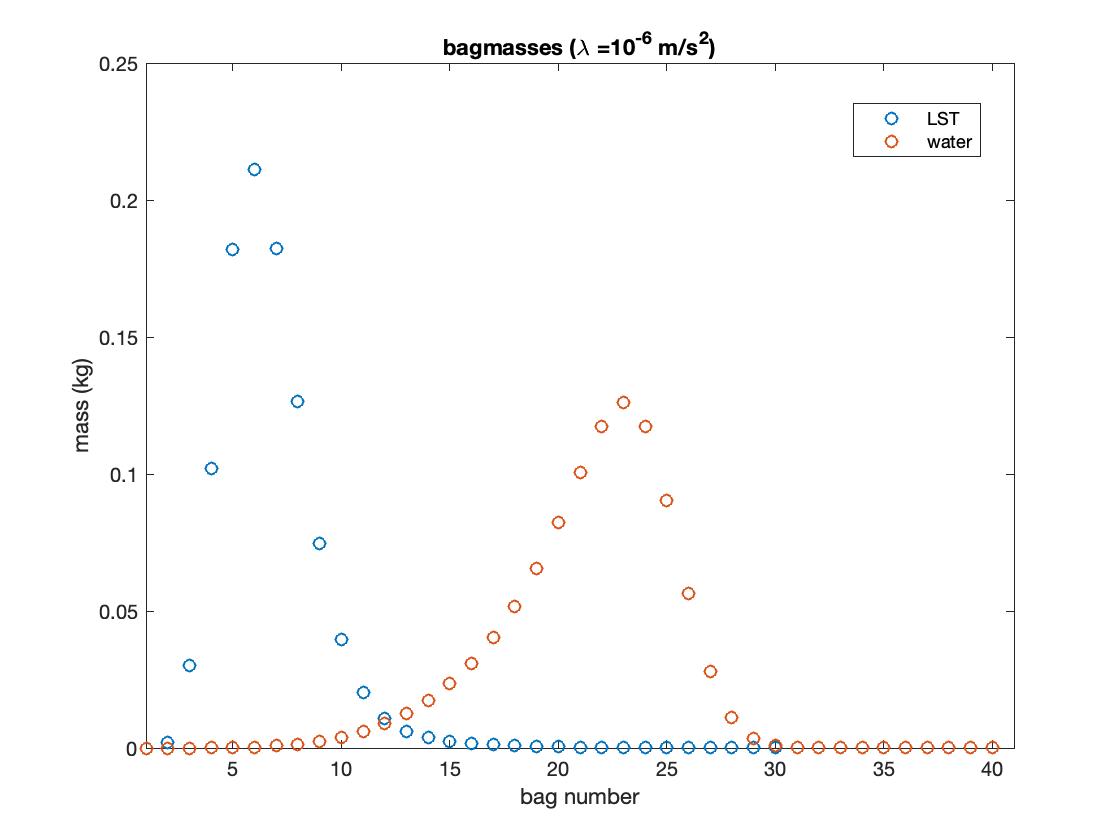}
\caption{Plots of the model bag masses $m_i^\varepsilon$ $(1\leq i\leq 40,\ \varepsilon\in\{1,2\})$ computed via (\ref{model masses}). The blue data points are those for LST while the red data points are those for water. Bags are taken at uniformly spaced times starting at $t=0$ and ending at $12.7$ hours.
}\label{fig: bag masses water-LST}
\end{center}
\end{figure}

Bag masses are modelled by the integrals (\ref{model masses}), the value of which are displayed in Figure \ref{fig: bag masses water-LST} (blue for water and red for the LST solution) with $\lambda =10^{-6}\,$m/s${}^2$. Here we see the effect of the more dense LST solution elutriating particles earlier than water due to its enhanced buoyancy effect.

\begin{figure}[h]
    \centering
    \subfloat{{\includegraphics[width=8cm]{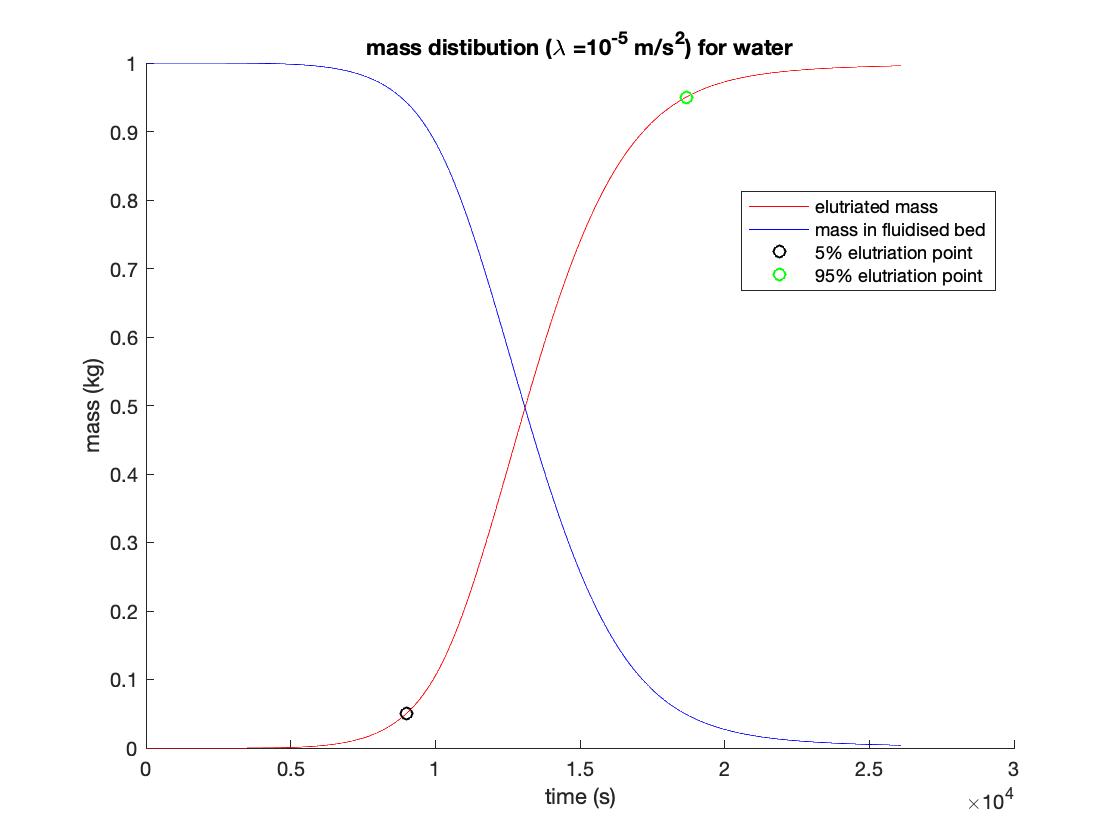} }}%
    \subfloat{{\includegraphics[width=8cm]{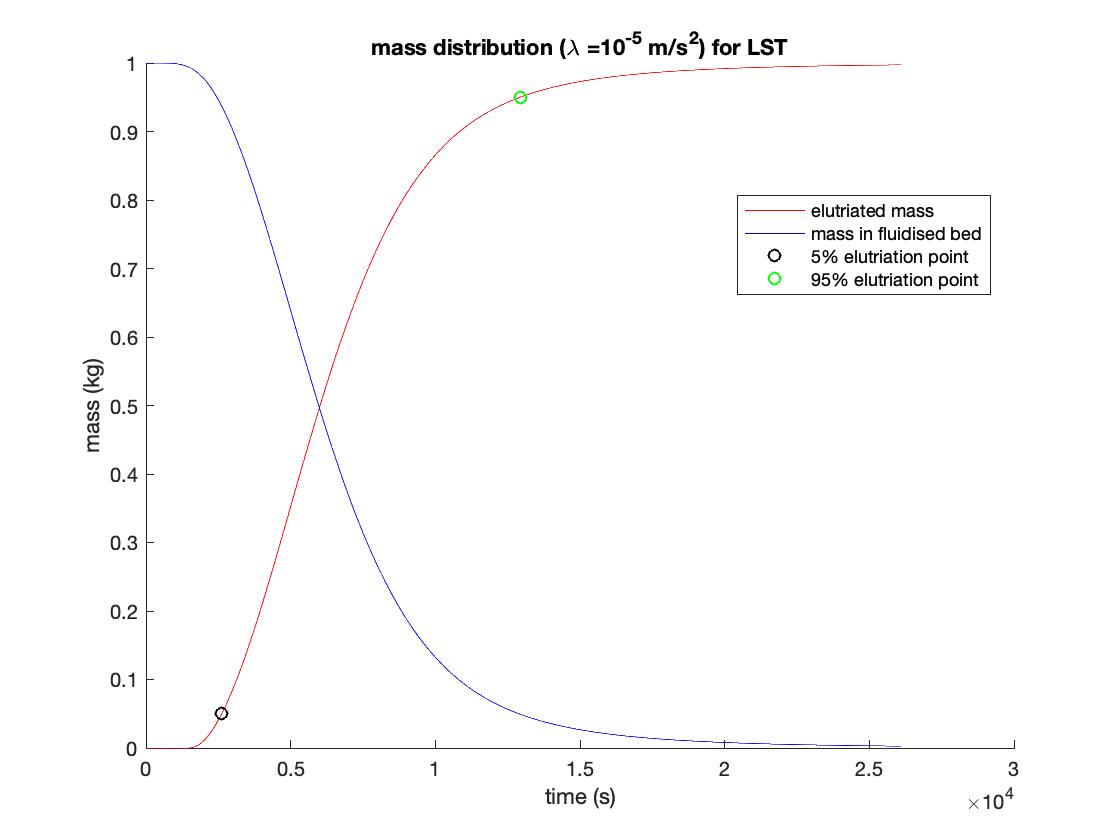}}}%
    \caption{The mass remaining in the fluidised bed at time $t$ (blue) and the total mass that has elutriated from the system by time $t$ (red) with $\lambda=10^{-5}\,$m/s${}^2$. The image on the left shows the results for the water experiment and the image on the right refers to the results for the LST experiment. Points at which $5\%$ $(T_0)$ of the initial mass has elutriated are coloured black and points where $95\%$ ($T_N$) of the initial mass has elutriated are coloured green.}%
    \label{fig: mass balance}%
\end{figure}

The mass balance between the fluidised bed, the channels and the elutriated material is demonstrated in Figure \ref{fig: mass balance}. The images show (in blue) the mass in the fluidised bed at time $t$ given by 
$$m_b(t)=\int_0^tM(s,0)e^{-kd_s(t)}\, ds$$
and (in red) the total mass that has elutriated by time $t$ given by
$$m_e(t)=\int_{d_s(t)>\ell} M(s,0)(1-e^{-k(d_s(t)-\ell )})\, ds$$
where $\ell$ is the length of the inclined channels. Here $\lambda=10^{-5}\,$ m/s${}^2$. The image on the left relates to the water experiment and the image on the right relates to the LST experiment. The $5\%$ and $95\%$ elutriation points are marked on each elutriation plot. The runtime of the water experiment (i.e., the difference between the time of $95\%$ elutriation and the time of $5\%$ elutriation) is $2.69$ hours, while the runtime of the LST experiment is $2.86$ hours.

The number $N$ of basis functions $\{\phi_j\}_{j=1}^N$ in (\ref{M exp}) should be as large as possible to obtain good fidelity of the solution and to cover the support of the mass distribution function. In real experiments, this number is determined by the method used in sizing the feed. On the other hand, there are restrictions on the number $p$ of elutriated sample bags that can be collected. The more samples collected (the higher the value of $p$), the smaller on average each one will be. This increases the noise to signal ratio in the raw data, as well as the amount of post-experiment sample handling and analysis that needs to be undertaken. The number of flow increments used in previously published elutriation work has varied from $6$ to $15$ (\cite{lowes et al 2020}, \cite{walton et al 2010}). Typically $N$ will be larger than  $p$. 

We restrict attention to the case where the $\phi_j$ are either splines of order $0$ (step functions) or $1$ (piecewise linear splines) for in these cases the non-negativity of $M$ as in (\ref{M exp}) is equivalent to the non-negativity of the coefficients $\{a_j\}_{j=1}^N$.

Let $s_1<s_2<\cdots <s_N<\infty$. The splines of order $0$ with respect to the {\it knots} $\{s_j\}_{j=1}^N$ are the piecewise constant functions
$$\phi_j^{(0)}(s)={\mathbf 1}_{[s_{j},s_{j+1})}(s)=\begin{cases}1&\text{ if $s_{j}\leq s<s_{j+1}$}\\
0&\text{ else.}
\end{cases}$$
The splines of order $1$ with respect to the knots $\{s_i\}_{i=1}^{N-1}$ are the piecewise linear functions given by
$$\phi^{(1)}_j(s)=\begin{cases}\dfrac{s-s_{j}}{s_{j+1}-s_{j}}&\text{ if $s_{j}\leq s\leq s_{j+1}$}\\
\dfrac{s_{j+2}-s}{s_{j+2}-s_{j+1}}&\text{ if $s_{j+1}\leq s\leq s_{j+2}$}\\
0&\text{ else.}
\end{cases}$$
Note that for $i=0,1$, the splines take non-negative values only and are {\it cardinal} in the sense that 
$$\phi^{(i)}_{j+1}(s_k)=\begin{cases}
1&\text{ if }j=k\\
0&\text{ else.}\end{cases}$$

In this paper, most of our computations are performed with respect to  splines of order zero -- which generate (discontinuous) step functions -- rather than the splines of order one -- which generate (continuous) piecewise linear functions. We do show the result of a deconvolution using linear splines to show that they allow for greater fidelity in the approximation to the continuous feed.

The placement of the knots are determined by the means by which we generate experimental size distributions -- either through the use of sieves, or through laser sizing. Both methods generate knots that are geometrically distributed, i.e., there exist $\beta >0$ and  $\delta >1$ such that $s_j=\beta\alpha^{j-1}s_1$ $(1\leq j\leq N)$. We conform to the values of $\alpha$ and $\beta$ given in Section \ref{subsec: size dist data}.

When using splines of order zero, the computation of the cross-grammian matrix entries $B_{ij}$ of (\ref{B entries}), the matrix entries $B_{ij}$ from (\ref{B entries}) simplify to
$$B_{ij}=\int_{s_j}^{s_{j+1}}L_i(s)\, ds$$
and the coefficients $\{c_j\}_{j=1}^N$ of (\ref{mass con}) become
$c_j=s_{j+1}-s_{j}$.

\subsection{Numerics of the inverse problem}

In this section we apply the theory developed in Section \ref{sec: deconvolution -- the inverse problem} to the inversion of the bag mass data obtained from (\ref{model masses}). The MATLAB function {\tt{lsqlin}} is used to solve the constrained optimisation problem (\ref{constrained opt problem}). We also incorporate {\it a posteriori} information about the feed -- namely that it takes the value $0$ at its largest scale -- into the problem, so that (\ref{constrained opt problem}) becomes
\begin{align}
&\text{minimise\ }\|B{\mathbf a}-{\mathbf m}_e\|^2+\alpha\|Q{\mathbf a}\|^2\notag\\
&\text{subject to }{\mathbf c}^{\text T}{\mathbf a}=\bar m\text{ and }a_j\geq 0\text{ for }1\leq j\leq N\text{ and }a_N=0.\label{constrained opt problem 2}
\end{align}
The extra constraint, coupled with the regularisation, promotes smooth behaviour of the solution at the end of its support, and reduces the space of solutions in which we search. Note that we make no assumption about the feed size distribution other than it is non-negative, the total mass is known, and we have a bound on the size of the largest particle. 

Reconstructions of the feed are computed using the value of the regularisation parameter that gives the lowest relative error as measured (as a percentage) by
$$\text{rel. err.}=\left[\frac{\int_0^z(M^{\#}(s)-\bar M(s))^2\, ds}{\int_0^zM^{\#}(s)^2\, ds}\right]^{1/2}\times 100$$ 
where $\bar M(s)=\sum_{j=1}^Na_j\phi_j(s)$ is the reconstructed feed arising from the optimisation, and $M^{\#}(s)=\sum_{j=1}^Nb_j\phi_j(s)$ is the spline that best approximates the feed in the least squares sense. In the case of splines of order zero, the coefficients are given by $b_j=\frac{1}{s_{j+1}-s_j}\int_{s_j}^{s_{j+1}}M(s)\, ds$.

\begin{figure}[h]
\begin{center}
\includegraphics[width=12cm,height=10cm]{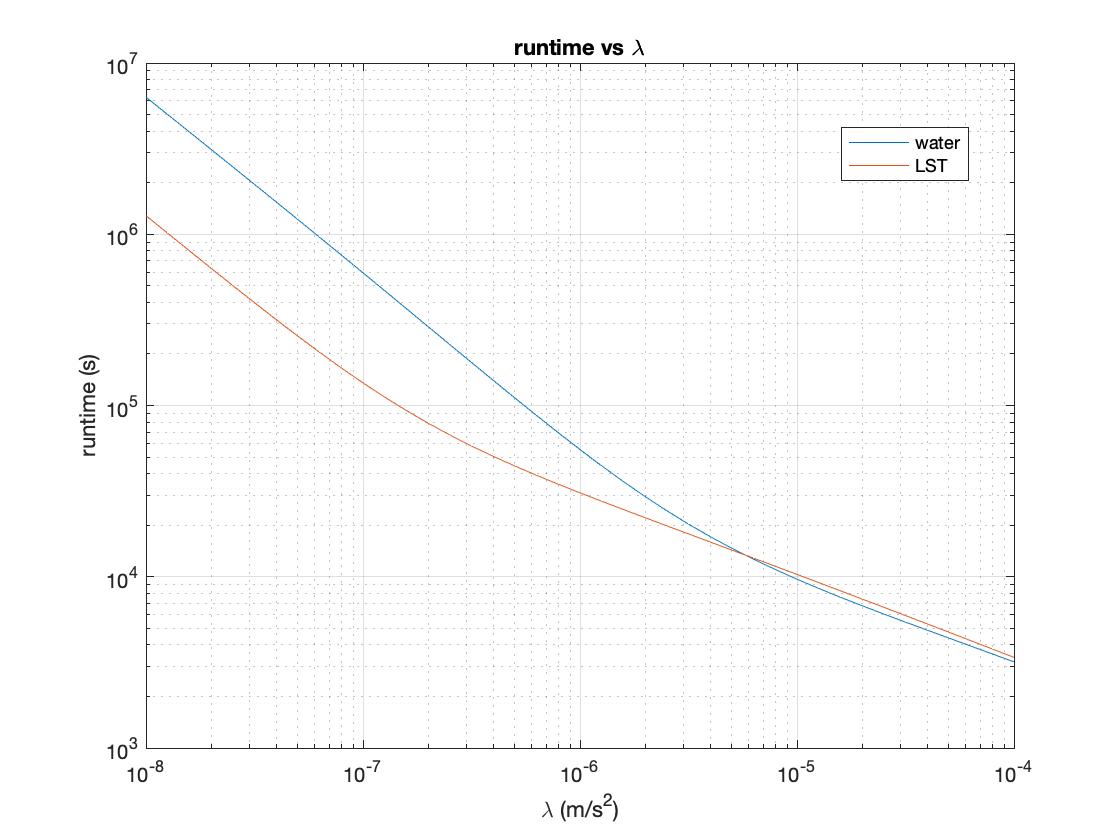}
\caption{Log-log plot of runtime against $\lambda $. High values of $\lambda$ give comparable runtimes for water and LST while for lower values of $\lambda$, the LST runtimes are significantly shorter.
}\label{fig: runtime vs lambda loglog}
\end{center}
\end{figure}

The effect of the fluid flow acceleration parameter $\lambda$ on runtime is very significant. Lower values of $\lambda$ produce better reconstructions at the expense of increased runtime. Figure \ref{fig: runtime vs lambda loglog} shows that for high values of $\lambda$, the runtimes of water and LST experiments are comparable (with water runs slightly shorter) while for lower values of $\lambda$ the runtime of LST experiments is much lower.

\begin{figure}%
    \centering
    \subfloat{{\includegraphics[width=8cm]{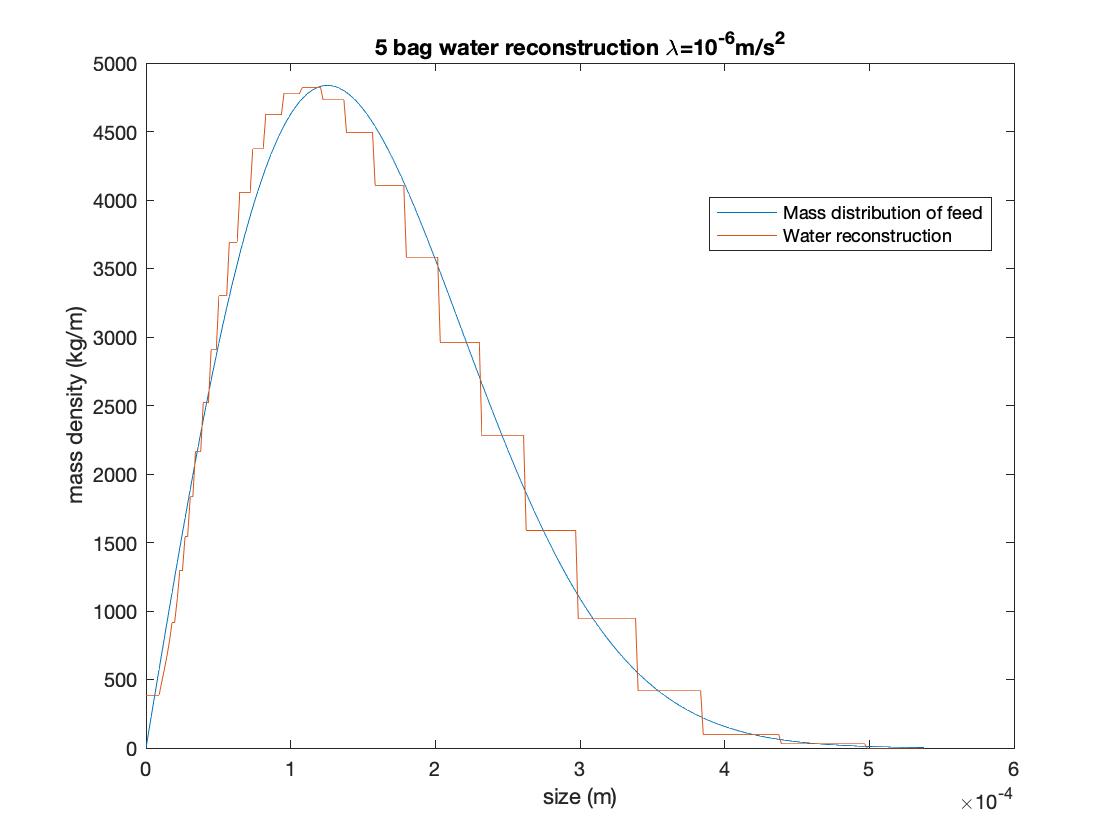} }}%
    \subfloat{{\includegraphics[width=8cm]{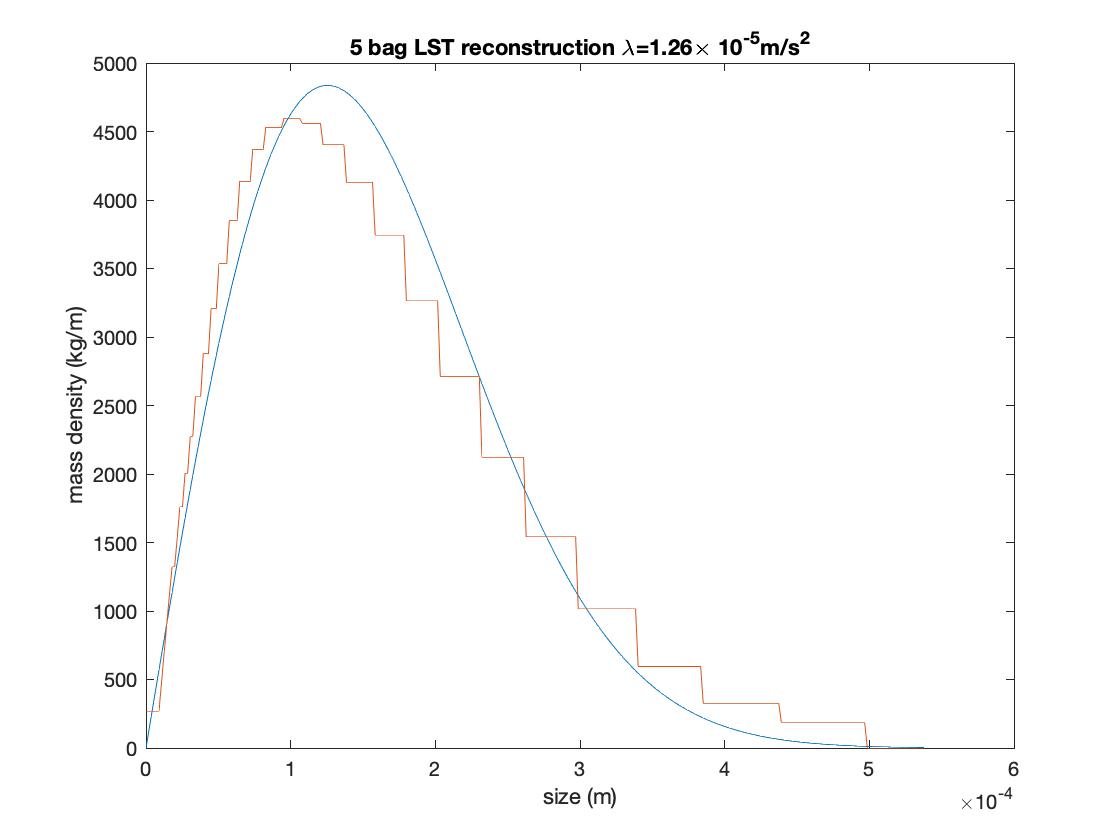} }}%
    \caption{Plot of the synthetic feed (blue) and the reconstruction of the feed (red) from $5$ bag experiments with runtime $2.69$ hours. The water reconstruction (relative error $5.4\%$) is on the left and the LST reconstruction (relative error $11.6\%$) is on the right. Best regularisation constants are used in both cases.}. %
    \label{fig: 5 bag reconstructions runtime =2.67 hours}%
\end{figure}

\begin{figure}%
    \centering
    \subfloat{{\includegraphics[width=8cm]{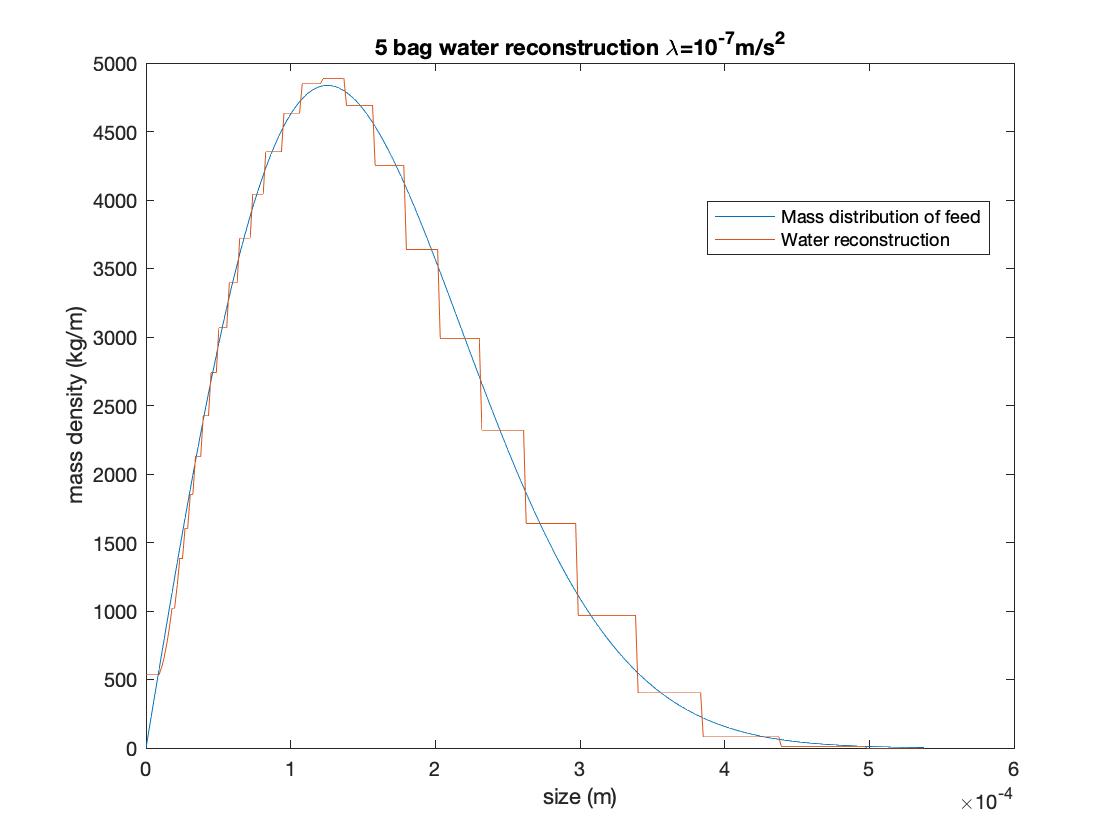} }}%
    \subfloat{{\includegraphics[width=8cm]{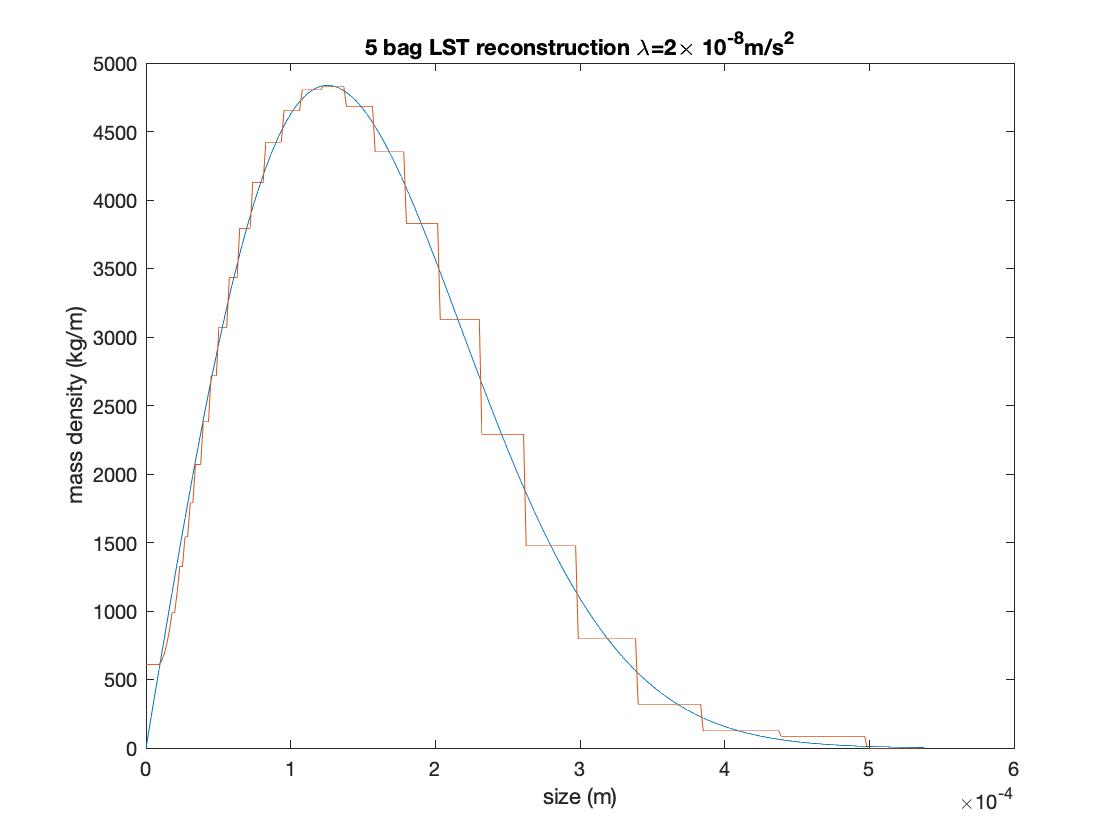} }}%
    \caption{Plot of the synthetic feed (blue) and the reconstruction of the feed (red) from $5$ bag experiments with runtime $165.2$ hours. The water reconstruction (relative error $3.5\%$) is on the left and the LST reconstruction (relative error $2.4\%$) is on the right. Best regularisation constants are used in both cases.}. %
    \label{fig: 5 bag reconstructions runtime =165.2 hours}%
\end{figure}

\begin{figure}%
    \centering
    \subfloat{{\includegraphics[width=8cm]{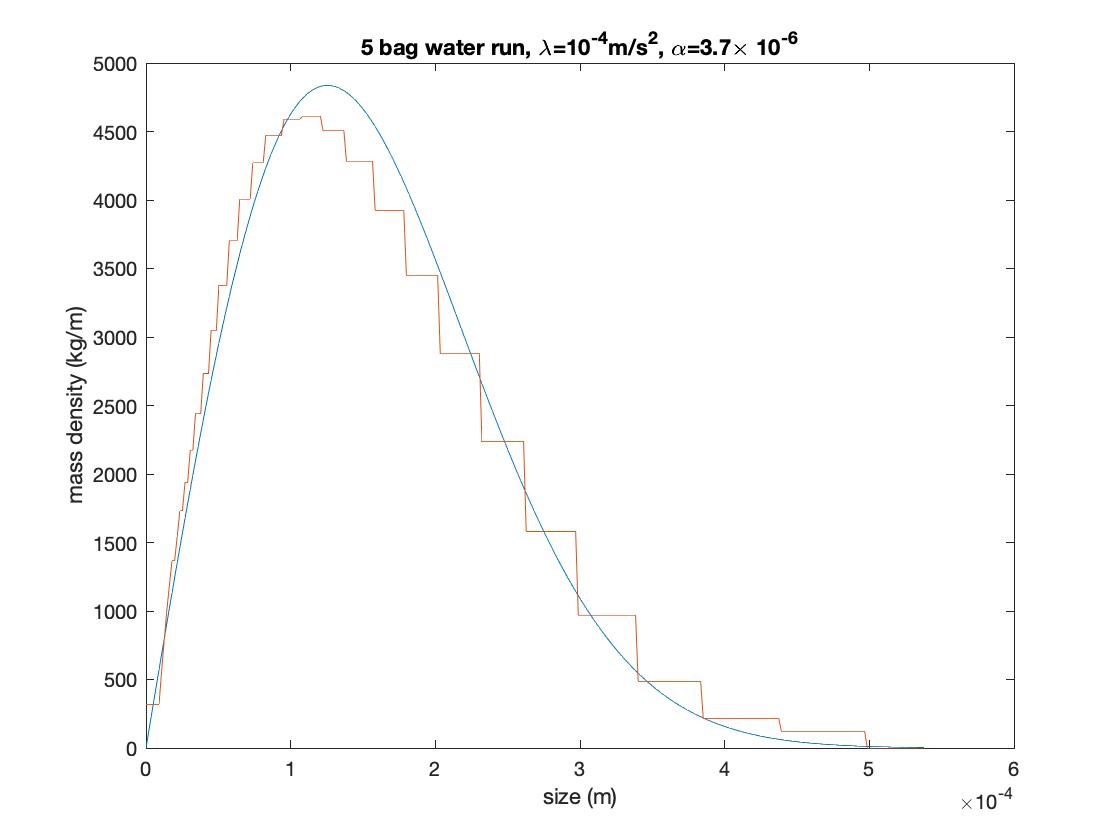} }}%
    \subfloat{{\includegraphics[width=8cm]{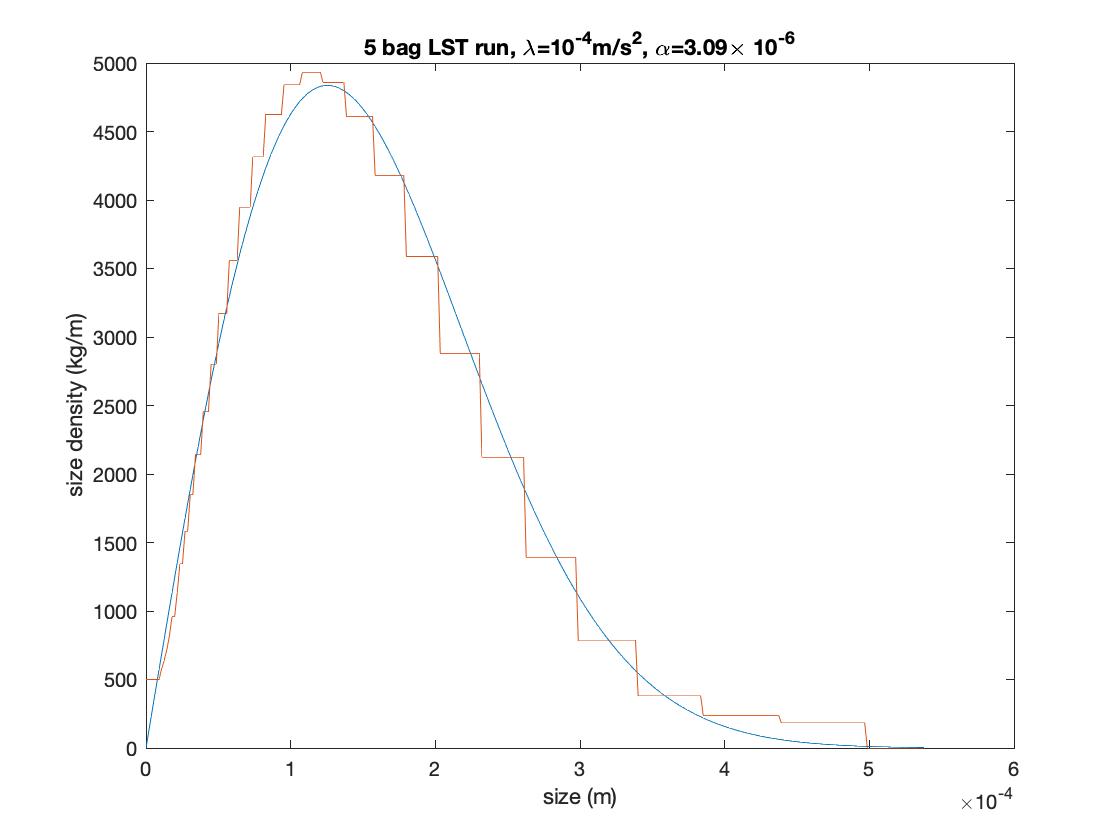} }}%
    \caption{Plot of the synthetic feed (blue) and the reconstruction of the feed (red) from $5$ bag experiments with $\lambda=10^{-4}$. The water reconstruction (relative error $7.8\%$) is on the left and the LST reconstruction (relative error $5.4\%$) is on the right. Best regularisation constants are used in both cases.}. %
    \label{fig: 5 bag reconstructions lambda=1e-4}%
\end{figure}

\begin{figure}[h]
\begin{center}
\includegraphics[width=12cm,height=10cm]{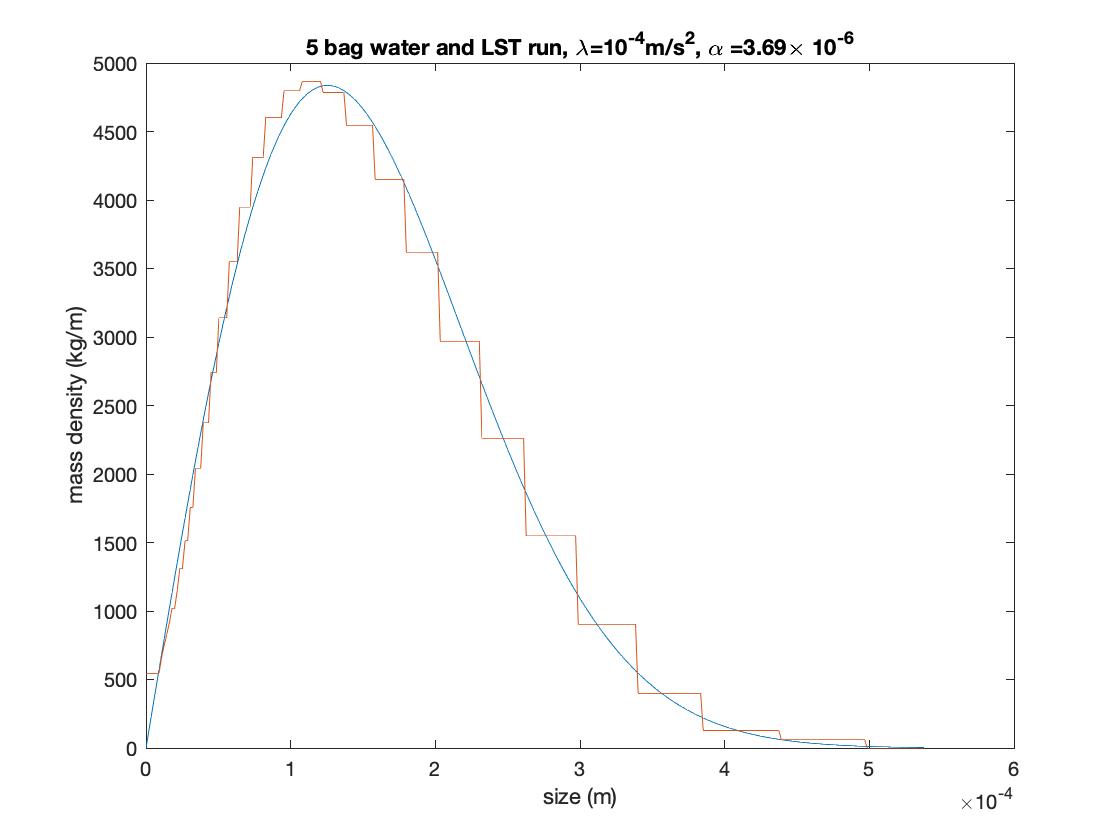}
\caption{Plot of the reconstruction obtained using the bag weight data generated in the $5$ bag water run and the $5$ bag LST run of the experiment with $\lambda=10^{-4}\,$m/s${}^2$ and best regularisation constant. The relative error is $4.4\%$.}\label{fig: 5 bag water+LST reconstructions lambda=1e-4}
\end{center}
\end{figure}

\begin{figure}%
    \centering
    \subfloat{{\includegraphics[width=8cm]{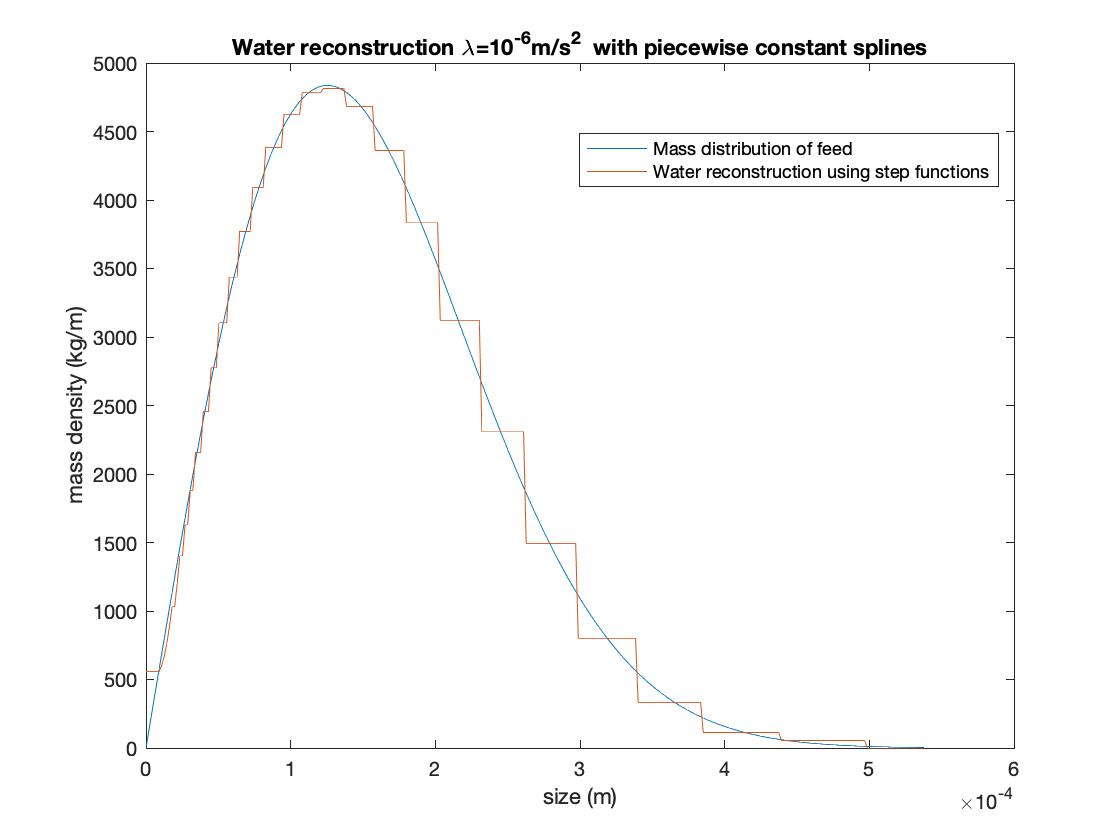} }}%
    \subfloat{{\includegraphics[width=8cm]{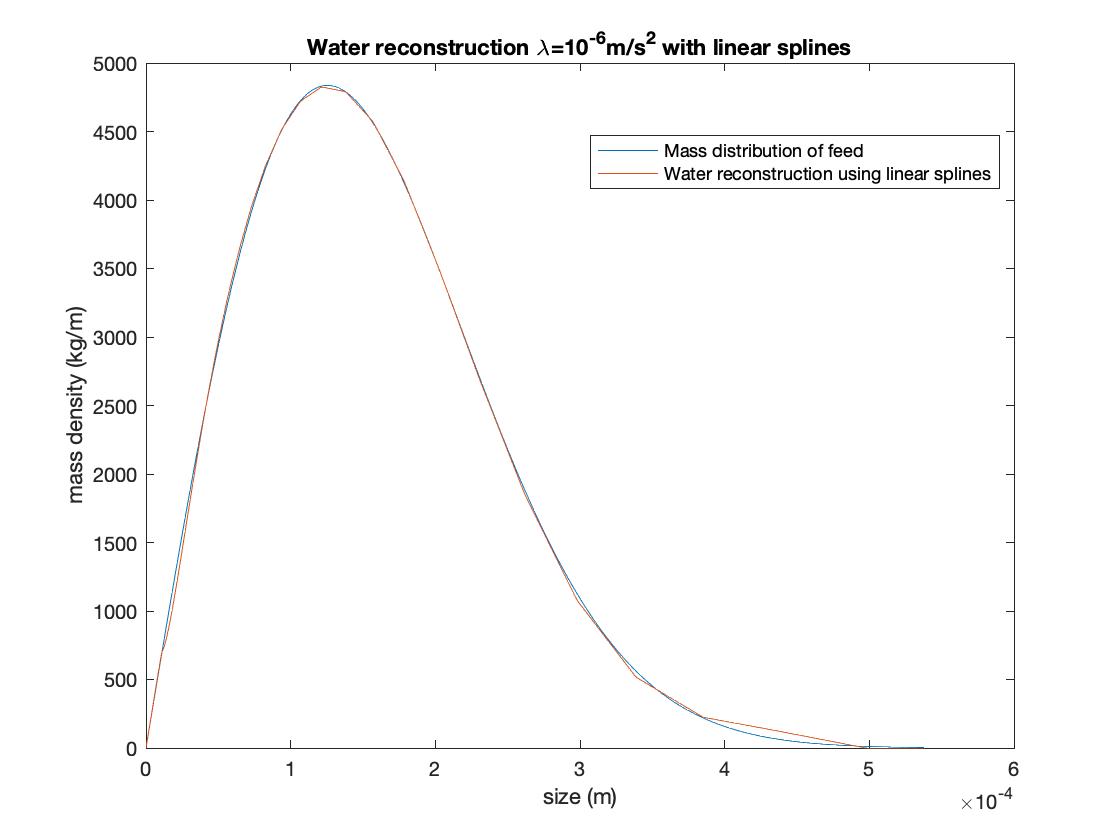} }}%
    \caption{Plots of water reconstructions with $\lambda =10^{-6}$ and best regularising constants using splines of order $0$ (left) and splines of order $1$ (right)}. %
    \label{fig: piecewise const v linear approx}%
\end{figure}

To assess the merits of using water versus LST, we present the following observations.
\begin{enumerate}

\item[(i)] In a water run of the experiment with $\lambda =10^{-5}\,$m/s${}^2$ we find that the runtime is $2.69$ hours. The best regularisation constant for these parameters is $\alpha =3.51\times 10^{-6}$ and the water reconstruction with these parameters is given in the image on the left in Figure \ref{fig: 5 bag reconstructions runtime =2.67 hours}. This reconstruction has relative error of $5.4\%$. To achieve the same runtime in an LST experiment, the value of $\lambda$ becomes $1.26\times 10^{-5}\,$m/s${}^2$. The best regularisation parameter is $\alpha =3.79\times 10^{-6}$ and the LST reconstruction with these parameters is given in the image on the right in Figure \ref{fig: 5 bag reconstructions runtime =2.67 hours}. This reconstruction has a relative error of $11.6\%$.

\item[(ii)] In a water run of the experiment with $\lambda =10^{-7}\,$m/s${}^2$ we find that the runtime is $165.2$ hours. The best regularisation constant for these parameters is $\alpha =6.07\times 10^{-6}$ and the water reconstruction with these parameters is given in the image on the left in Figure \ref{fig: 5 bag reconstructions runtime =165.2 hours}. This reconstruction has relative error of $3.5\%$. To achieve the same runtime in an LST experiment, the value of $\lambda$ becomes $1.2\times 10^{-8}\,$m/s${}^2$. The best regularisation parameter is $\alpha =5.03\times 10^{-6}$ and the LST reconstruction with these parameters is given in the image on the right in Figure \ref{fig: 5 bag reconstructions runtime =165.2 hours}. This reconstruction has a relative error of $2.4\%$.

\item[(iii)] We investigate the benefits of combining the data from water and LST runs of the experiment. With $\lambda=10^{-4}\,$m/s${}^2$, Figure \ref{fig: 5 bag reconstructions lambda=1e-4} shows reconstructions from a water run (left) with $\alpha =3.7\times 10^{-6}$ and relative error $7.8\%$ and an LST run (right) with $\alpha =3.09\times 10^{-6}$ and relative error $5.5\%$. In Figure \ref{fig: 5 bag water+LST reconstructions lambda=1e-4} we see a reconstruction from the combined data generated by the water and LST runs, now with $\alpha =3.69\times 10^{-6}$ and a relative error of $4.4\%$. Note the improved performance of this calculation for large particle sizes.

\item[(iv)] As described in Section \ref{subsec: forward prob}, rather than using step functions (splines of order $0$) as basis functions in the model (\ref{M exp}) for the feed we may use piecewise linear  functions (splines of order $1$). The advantage of the splines of order $1$ lies in their smoothness which allows for a better least squares fit to the assumed continuous feed. In Figure \ref{fig: piecewise const v linear approx} we show the reconstruction of the feed for a water run of the experiment with $\lambda=10^{-6}\,$m/s${}^2$ and best regularisation parameters using splines of order $0$ (left, with relative error $1.8\%$) and splines of order $1$ (right with relative error $1.3\%$ ). The piecewise linear approximation matches the continuous feed much more closely than the piecewise constant approximation since (unlike the best piecewise constant spline approximation of the feed) the best piecewise linear spline approximation of the feed is almost indistinguishable from the continuous feed.
\end{enumerate}

Comparing Figures \ref{fig: 5 bag reconstructions runtime =2.67 hours} and \ref{fig: 5 bag reconstructions runtime =165.2 hours}, we see that for both water and LST increasing the runtime (lowering the value of $\lambda$) improves the accuracy of deconvolution. This is due to the effect shown in Figure \ref{fig: kernel plots} -- lowering the value of $\lambda$ causes the range of particle sizes in the bags to narrow, thus giving greater resolution and accuracy of deconvolution. There is therefore a trade-off between accuracy and runtime.

\subsection{Noisy measurements}
\begin{figure}%
    \centering
    \subfloat{{\includegraphics[width=8cm]{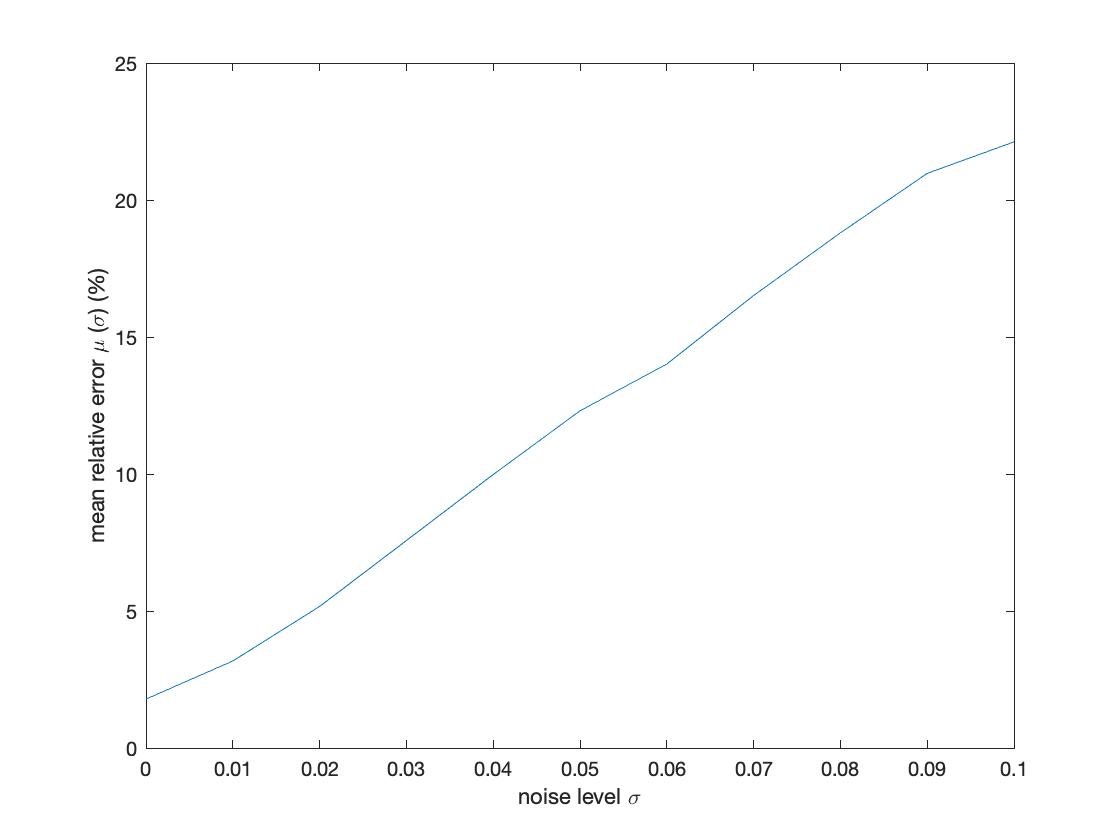} }}%
    \subfloat{{\includegraphics[width=8cm]{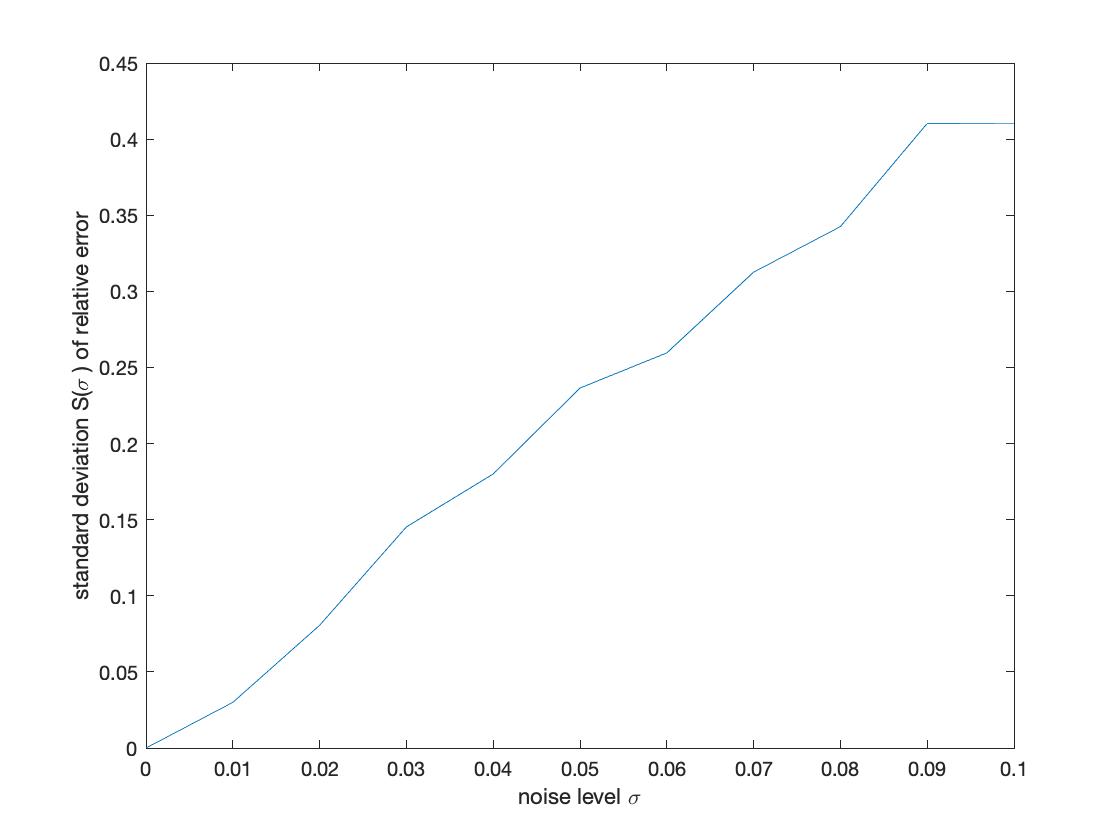} }}%
    \caption{Plots of the the sample mean error $\mu$ and standard deviation $S$ against $\sigma$ as defined in (\ref{mean and std}) with $500$ simulations.}%
    \label{fig: effect of noise}%
\end{figure}

As a way of modelling real data, we add Gaussian noise to the computed bag masses and perform the reconstruction on the noisy data. In particular, given $\sigma >0$ and with $m_{e,i}$ as in (\ref{model masses}), we let
$$m_{e,i}^{\sigma}=\max (m_{e,i} +N(0,\sigma m_{e,i}),0)$$
where $N(0,\sigma m_{e,i})$ is a random number generated by the normal distribution with mean $0$ and standard deviation $\sigma m_{e,i}$. Then we compute $N$ reconstructions $\bar M_n^\sigma (s)$ $(1\leq n\leq N)$ and the errors $E_n(\sigma )=(\int_0^\infty (M(s)-\bar M_n^\sigma (s))^2\, ds)^{1/2}$. We report on the mean error $\mu (\sigma )$ and the standard deviation $S(\sigma )$ given by 
\begin{equation}
\mu (\sigma )=\frac{1}{N}\sum_{i=1}^{N}E_n(\sigma );\quad  S(\sigma ) =\frac{1}{N}(\sum_{i=1}^{N}(E_n(\sigma )-\mu (\sigma ))^2)^{1/2}.\label{mean and std}
\end{equation}
These quantities are plotted in Figure \ref{fig: effect of noise} for $N=500$. It appears that (at least for small values of $\sigma$), the mean error and standard deviation increase linearly with $\sigma$.

\section{Conclusion}
This work has shown that bag mass data from the batch fractionation of a mono-density feed sample can be deconvoluted to infer the feed particle size distribution, based on a known model for the system’s behaviour. Whilst the model used in this work will not be a perfect representation of a real system, it nevertheless does provide some practical insights into how these systems function. It was shown that using slower flow acceleration rates gives improved resolution and also that the error in the result only grows linearly with errors in the measured bag weights. Hence this means that one can hope to achieve whatever level of accuracy is desired, simply by running the experiments more slowly. It was also shown that combining test data from two different liquids may also improve the accuracy. The benefits of this latter approach are likely to become more significant when considering feed samples which vary in both size and density, which is the focus of ongoing work.


\begin{thebibliography}{10}

\bibitem{amariei et al} D. Amariei, D. Michaud, G. Paquet, G. M. Lindsay. \newblock The use of the Reflux Classifier for iron ores: assessment of fine particles recovery at pilot scale,\newblock {Min. Eng.} (2014) 62, 66–73 {\tt https://doi.org/10.1016/j.mineng.2013.11.011}

\bibitem{amoash et al 2025} M.E. Amosah, J. Zhou, K.P. Galvin. \newblock Fourth generation gravity separation using the REFLUX™ Classifier,\newblock {Min. Eng.} (2025)  224 109216 {\tt https://doi.org/10.1016/j.mineng.2025.109216}

\bibitem{bird et al 1976} B. Bird, W.E. Stewart, E.N. Lightfoot. \newblock Transport Phenomena. New York: Wiley (1976).

\bibitem{boycott}
A.E. Boycott.
\newblock Sedimentation of blood corpuscles,
\newblock {Nature} (1920) 104, 532 {\tt https://doi.org/10.1038/104532b0}

\bibitem{CCC} Central Coast Consulting.
\newblock LST Heavy Liquid – Comparison with other heavy liquids,
\newblock Available online: {\tt https://www.chem.com.au/comparison.html} (accessed on 15 Dec 2023)
	
	
	\bibitem{chu et al 2020} H. Chu, Y. Wang, D. Lu, Z. Liu, X. Zheng. \newblock Pre-concentration of fine antimony oxide tailings using an agitated reflux classifier,\newblock {Powder Technol.} (2020), 376, 565-572 {\tt https://doi.org/10.1016/j.powtec.2020.08.055}
	
	\bibitem{fox and mcdonald 1998} R.W. Fox, A.T. McDonald. \newblock Introduction to Fluid Mechanics (1998) (5th ed.). John Wiley \& Sons
	
	\bibitem{galvin et al 2018} K.P. Galvin, S.M. Iveson, D.M. Hunter. \newblock Deconvolution of fractionation data to deduce consistent washability and partition curves for a mineral,\newblock {Min. Eng.} (2018) 125 94-110 {\tt https://doi.org/10.1016/j.mineng.2018.05.010}
	
		\bibitem{galvin et al 2020b} K.P. Galvin, S.M.  Iveson, J. Zhou, C.P. Lowes. \newblock Influence of inclined channel spacing on dense mineral partition in a REFLUX™ Classifier. Part 1: Continuous steady state, {Min. Eng.} (2020) 146, 106112 {\tt https://doi.org/10.1016/j.mineng.2019.106112}
		
		\bibitem{galvin et al 2020c} K.P. Galvin, S.M. Iveson, J. Zhou, C.P. Lowes. \newblock Influence of inclined channel spacing on dense mineral partition in a REFLUX™ Classifier. Part 2: Water based fractionation, {Min. Eng.} (2020) 155, 106442 {\tt https://doi.org/10.1016/j.mineng.2020.106442}
	
	\bibitem{galvin-liu 2011} K.P. Galvin, H. Liu. \newblock Role of inertial lift in elutriating particles according to their density,\newblock {Chem. Eng. Sci.}, (2011) 66, 3686–3691. {\tt https://doi.org/10.1016/j.ces.2011.05.002}
	
	\bibitem{galvin et al 2009}
	\newblock K.P. Galvin, K. Walton,  J. Zhou. How to elutriate particles, according to their density,\newblock {Chem. Eng. Sci.} (2009) 64, 2003–2010. {\tt doi:10.1016/j.ces.2009.01.031}
	
	\bibitem{galvin et al 2016} K.P. Galvin, J. Zhou, A.J. Price, P. Agrwal, S.M. Iveson. \newblock Single-stage recovery and concentration of mineral sands using a REFLUX™ Classifier,\newblock {Min. Eng.} (2016), 93, 32-40. {\tt https://doi.org/10.1016/j.mineng.2016.04.010}

	
	\bibitem{galvin et al 2020a} K.P. Galvin, J. Zhou, J.L. Sutherland, S.M. Iveson.  \newblock Enhanced recovery of zircon using a REFLUX™ Classifier with an inclined channel spacing of 3 mm, {Min. Eng.} (2020) 147, 106148 {\tt https://doi.org/10.1016/j.mineng.2019.106148}
	
	
\bibitem{haidler-levenspiel 1989} A. Haider,  O. Levenspiel. \newblock Drag coefficient and terminal velocity of spherical and non-spherical particles, \newblock {Powder Technol.} (1989) 58, 63-70.
	
	\bibitem{hunter et al 2014} D.M. Hunter, S.M. Iveson, K.P. Galvin. \newblock The role of viscosity in the density fractionation of particles in a laboratory-scale REFLUX™ Classifier, \newblock {Fuel} (2014) 129, 188-196. {\tt http://dx.doi.org/10.1016/j.fuel.2014.03.063}
	
	\bibitem{iveson et al 2024} S.M. Iveson, N. Boonzaier, K.P. Galvin. \newblock Beneficiation of high-density Tantalum ore in the REFLUX(TM) concentrating classifier analysed using batch fractionation assay and density data,\newblock {Minerals} (2024) 14, 197 {\tt https://doi.org/10.3390/min14020197}
	
	\bibitem{lowes et al 2020}
	C. Lowes, J. Zhou and K.P. Galvin.
	\newblock Improved density fractionation of minerals in the REFLUX${}^\text{TM}$ Classifier using LST as a novel fluidising medium.
	\newblock {Min. Eng.} (2020) 146:106-145.
	
	\bibitem{mueller-siltanen}
	J.L. Mueller and S. Siltanen.
	\newblock Linear and Nonlinear Inverse Problems with Applications,
	\newblock (2012) SIAM, Philadelphia.
	
	\bibitem{rosin-rammler}
	P.Rosin and E.Rammler.
	\newblock The laws governing the fineness of powdered coal,
	\newblock {J. Inst. Fuel} (1993) 7: 29-36.
	
	
	\bibitem{starrett}
	J. Starrett.
	\newblock Density Fractionation of Mineral Ores Using Heavy Liquid Fractionation Through Inclined Channels,
	\newblock Honours thesis, (2021) University of Newcastle (Australia).
	
	\bibitem{strauss} W.A. Strauss.\newblock Partial Differential Equations: an Introduction, \newblock (2008) Wiley, Hoboken.
	
	\bibitem{walton et al 2010} K. Walton, J. Zhou, K.P. Galvin. Processing of fine particles using closely spaced inclined channels, \newblock {Adv. Powder Technology} (2010) 21, 386-391 {\tt http://dx.doi.org/10.1016/j.apt.2010.02.015}
	
	\bibitem{Wills}
	B.A. Wills.
	\newblock Mineral Processing Technology: An introduction to the practical aspects of ore treatment and mineral recovery (6th ed.)
	\newblock (1997) Butterworth-Heinemann.
	
	\bibitem{wright}
	B. Wright.
	\newblock Density Fractionation of a Binary Feed Through the Use of Inclined Channels,
	\newblock Honours thesis, (2021) University of Newcastle (Australia).
	
	\bibitem{zigrang-sylvester 1981}
	D.J. Zigrang and N.D. Sylvester.
	\newblock An explicit equation for particle settling velocities in solid-liquid systems,
	\newblock {AIChE Journal} (1081) 27 (6), 1043-1044 {\tt  https://doi.org/10.1002/aic.690270629}
	
\end{thebibliography}
\end{document}